# LARGE DEVIATIONS AND A KRAMERS' TYPE LAW FOR SELF-STABILIZING DIFFUSIONS

By Samuel Herrmann, Peter Imkeller and Dierk Peithmann

*Ecole des Mines de Nancy and Humboldt Universität zu Berlin*

We investigate exit times from domains of attraction for the motion of a self-stabilized particle traveling in a geometric (potential type) landscape and perturbed by Brownian noise of small amplitude. Self-stabilization is the effect of including an ensemble-average attraction in addition to the usual state-dependent drift, where the particle is supposed to be suspended in a large population of identical ones. A Kramers' type law for the particle's exit from the potential's domains of attraction and a large deviations principle for the self-stabilizing diffusion are proved. It turns out that the exit law for the self-stabilizing diffusion coincides with the exit law of a potential diffusion without self-stabilization and a drift component perturbed by average attraction. We show that self-stabilization may substantially delay the exit from domains of attraction, and that the exit location may be completely different.

## 1. Introduction.

We examine the motion of a particle subject to three sources of forcing. First, it wanders in a landscape whose geometry is determined by a potential. Second, its trajectories are perturbed by Brownian noise of a small amplitude. The third source of forcing can be thought of as self-stabilization. Roughly, it characterizes the influence of a large population of identical particles subject to the same laws of motion. They act on the individual through an attractive potential averaged over the whole population, which adds to the underlying potential drift. More formally, denote by $X_t^\varepsilon$ the random position of the particle at time $t$. It is governed by the $d$-dimensional SDE

$$(1.1) \qquad dX_t^\varepsilon = V(X_t^\varepsilon) \, dt - \int_{\mathbb{R}^d} \Phi(X_t^\varepsilon - x) \, du_t^\varepsilon(x) \, dt + \sqrt{\varepsilon} \, dW_t.$$

In this equation, $V$ denotes a vector field on $\mathbb{R}^d$, which we think of as representing a potential gradient, the first source of forcing. Without the other two sources the motion of the particle would just amount to the dynamical system given by the ODE

$$\dot{x} = V(x). \tag{1.2}$$

The small stochastic perturbation by Brownian noise $W$ of intensity $\varepsilon$ accounts for the second source of forcing. It is responsible for random behavior of $X^\varepsilon$, and allows for transitions between otherwise energetically unreachable domains of attraction. The integral term involving the process's own law $u_t^\varepsilon$ introduces a feature that we call *self-stabilization*. The distance between the particle's instantaneous position $X_t^\varepsilon$ and a fixed point $x$ in state space is weighed by means of a so-called *interaction function* $\Phi$ and integrated in $x$ against the law of $X_t^\varepsilon$ itself. This effective additional drift can be seen as a measure for the average attractive force exerted on the particle by an independent copy of itself through the attraction potential $\Phi$. In effect, this forcing makes the diffusion inertial and stabilizes its motion in certain regions of the state space.

Equations of the type (1.1) are obtained as mesoscopic limits of microsystems of interacting particles, as the number of particles in an ensemble of identical ones tends to infinity, and subject to the same first two sources of forcing, that is, the force field $V$ and the Brownian noise of intensity $\varepsilon$. Suppose we are given an interaction function $\Phi$, that is, for any two particles located at $x$ and $y$ in state space the value $\Phi(x - y)$ expresses the force of mutual attraction. This attraction can, for instance, be thought of as being generated by electromagnetic effects. The dynamics of a particle system consisting of $N$ such particles is described by the stochastic differential equation

$$dX_t^{i,N} = V(X_t^{i,N})\,dt - \frac{1}{N}\sum_{j=1}^{N}\Phi(X_t^{i,N} - X_t^{j,N})\,dt + \sqrt{\varepsilon}\,dW_t^i,$$
$$X_0^{i,N} = x_0^i, \qquad 1 \le i \le N. \tag{1.3}$$

Here the $W^i$ are independent Brownian motions. The self-stabilizing effect we are interested in originates in the global action of the system on the individual particle motion in the large particle limit $N \to \infty$. Under suitable assumptions, in this limit the empirical measures $\frac{1}{N}\sum_{j=1}^{N}\delta_{X_t^{j,N}}$ can be shown to converge to some law $u_t^\varepsilon$ for each fixed time and noise intensity, and each individual particle's motion converges in probability to the solution of the diffusion equation

$$dX_t^i = V(X_t^i)\,dt - \int_{\mathbb{R}^d}\Phi(X_t^i - x)\,du_t^\varepsilon(x)\,dt + \sqrt{\varepsilon}\,dW_t^i. \tag{1.4}$$

The aim of this paper is to extend the well-known Kramers–Eyring law of exit from domains with noncritical boundaries by particles diffusing in potential landscapes with small Gaussian noise to systems (1.1) which include the described self-stabilization effect. In the potential gradient case without interaction, in which the individual particle's motion is interpreted by the solution trajectories $Z^\varepsilon$ of the SDE

$$(1.5) \qquad dZ_t^\varepsilon = -\nabla U(Z_t^\varepsilon)\, dt + \sqrt{\varepsilon}\, dW_t,$$

Kramers' law states that, in the small noise limit $\varepsilon \to 0$, the asymptotic exit time of $Z^\varepsilon$ from a potential well of height $H$ is of the order $\exp\{\frac{2H}{\varepsilon}\}$. See the beginning of Section 4 for a precise formulation of this. We derive a similar statement for self-stabilizing diffusions. In particular we examine how self-stabilization adds inertia to the individual particle's motion, delaying exit times from domains of attraction and altering exit locations. Mathematically, the natural framework for such an analysis is large deviations theory for diffusions. Our key ingredient for an understanding of the small noise asymptotics of the exit times proves to be a large deviations principle for self-stabilizing diffusions (1.1). In the potential gradient case, the rate function in the large deviations principle just minimizes the energy needed to travel in the potential landscape. If the particle undergoes self-stabilization, energy has to be minimized in a landscape which additionally takes into account the potential of an attractive force that depends on the particle's distance from the corresponding deterministic path (1.2). Our main results (Theorems 3.4 and 4.2, 4.3) state that the large deviations and the exit behavior of $X^\varepsilon$ are governed by this modified rate function. The techniques we employ to relate this time-inhomogeneous case to the classical time-homogeneous one stipulate the assumption that the boundaries of the domains avoid critical points of the potential.

Interacting particle systems such as (1.3) have been studied from various points of view. A survey about the general setting for interaction (under global Lipschitz and boundedness assumptions) may be found in [14]. There the convergence of the particle system to a self-stabilizing diffusion is described in the sense of a McKean–Vlasov limit, and asymptotic independence of the particles, known under the name *propagation of chaos*, as well as the link to Burgers' equation are established. Large deviations of the particle system from the McKean–Vlasov limit were investigated by Dawson and Gärtner [4]. Further results about the convergence of the empirical distribution of the particle system to the law of the self-stabilizing diffusion may be found in [3] or [10].

McKean studies a class of Markov processes that contains the solution of the limiting equation under global Lipschitz assumptions on the structure of the interaction [11]. A strictly local form of interaction was investigated by Stroock and Varadhan in simplifying its functional description to a Dirac measure [13]. Oelschläger studies the particular case where interaction is represented by the derivative of the Dirac measure at zero [12]. Funaki addresses existence and uniqueness for the martingale problem associated with self-stabilizing diffusions [9].

The behavior of self-stabilizing diffusions, in particular the convergence to invariant measures, was studied by various authors under different assumptions on the structure of the interaction; see, for example, [1, 2, 15, 16].

The material in this paper is organized as follows. In Section 2 we discuss existence and uniqueness of strong solutions to (1.1). Strong solvability is nontrivial in our setting due to the self-stabilizing term, and is required for the subsequent investigation of large deviations. In Section 3 we derive and analyze the rate function modified by self-interaction, and this way obtain a large deviations principle for the diffusion (1.1). This proves to be the key ingredient for the analysis of exit times and a derivation of a version of Kramers' law for self-stabilizing diffusions in Section 4. We conclude with an illustration of our main results by discussing some examples which emphasize the influence of self-stabilization on exit time and exit location (Section 5).

## 2. Existence and uniqueness of a strong solution.

The derivation of a large deviations principle for the self-stabilizing diffusion (1.1) in the subsequent section involves pathwise comparisons between diffusions in order to apply the usual tools from large deviations theory, such as contraction principles and the concept of exponential equivalence. Their applicability relies on strong existence and uniqueness for (1.1), which is nontrivial in our situation since the solution process's own law appears in the equation. The interesting interaction term $\int \Phi(X_t^\varepsilon - x) \, du_t^\varepsilon(x)$ also adds a considerable amount of complexity to the mathematical treatment. It depends on $u_t^\varepsilon = \mathbb{P} \circ (X_t^\varepsilon)^{-1}$; thus classical existence and uniqueness results on SDE as well as the classical results on large deviations for diffusions (Freidlin–Wentzell theory) are not directly applicable. Consequently, the question of existence and uniqueness of solutions for (1.1) is an integral part in any discussion of the self-stabilizing diffusion's behavior, and will be addressed in this section.

We follow Benachour et al. [1] to design a recursive procedure in order to prove the existence of the *interaction drift* $b(t, x) = \int \Phi(x - y) \, du_t^\varepsilon(y)$, the second drift component of (1.1). More precisely, we shall construct a locally Lipschitz drift term $b(t, x)$ such that the classical SDE

$$(2.1) \qquad dX_t^\varepsilon = V(X_t^\varepsilon) \, dt - b(t, X_t^\varepsilon) \, dt + \sqrt{\varepsilon} \, dW_t, \qquad t \geq 0,$$

admits a unique strong solution, which satisfies the additional condition

$$(2.2) \qquad b(t, x) = \int_{\mathbb{R}^d} \Phi(x - y) \, du_t^\varepsilon(y) = \mathbb{E}\{\Phi(x - X_t^\varepsilon)\}.$$

In (2.1) $W$ is a standard $d$-dimensional Brownian motion, and $V : \mathbb{R}^d \to \mathbb{R}^d$ mimics the geometrical structure of a potential gradient. Of course, for (2.1) to make sense, the drift $b$ has to be well defined, that is, the integral of (2.2) has to be finite, which depends upon certain moment conditions for $X^\varepsilon$ to be made precise later on. Apart from these moment conditions, existence and uniqueness for (1.1) will be understood in the sense that (2.1) and (2.2) hold with a unique $b$ and a pathwise unique process $X^\varepsilon$.

For locally Lipschitz interaction functions of at most polynomial growth, Benachour et al. [1] have proved the existence of strong solutions in the one-dimensional situation, and in the absence of the vector field $V$. Since $V$ forces the diffusion to spend even more time in bounded sets due to its dissipativity formulated below, it imposes no complications concerning questions of existence and uniqueness. Our arguments rely on a modification of their construction.

Besides some Lipschitz type regularity conditions on the coefficients, we make assumptions concerning the geometry of $V$ and $\Phi$ which render the system (3.1) dissipative in a suitable sense. All necessary conditions are summarized in the following assumption.

ASSUMPTION 2.1.

(i) The coefficients $V$ and $\Phi$ are locally Lipschitz, that is, for each $R > 0$ there exists $K_R > 0$ s.t.

$$(2.3) \qquad \|V(x) - V(y)\| + \|\Phi(x) - \Phi(y)\| \le K_R \|x - y\|$$

for $x, y \in B_R(0) = \{z \in \mathbb{R}^d : \|z\| < R\}$.

(ii) The interaction function $\Phi$ is rotationally invariant, that is, there exists an increasing function $\phi : [0, \infty) \to [0, \infty)$ with $\phi(0) = 0$ such that

$$(2.4) \qquad \Phi(x) = \frac{x}{\|x\|} \phi(\|x\|), \qquad x \in \mathbb{R}^d.$$

(iii) $\Phi$ grows at most polynomially: there exist $K > 0$ and $r \in \mathbb{N}$ such that

$$(2.5) \qquad \|\Phi(x) - \Phi(y)\| \le \|x - y\|(K + \|x\|^r + \|y\|^r), \qquad x, y \in \mathbb{R}^d.$$

(iv) $V$ is continuously differentiable. Let $DV(x)$ denote the Jacobian of $V$. We assume that there exist $K_V > 0$ and $R_0 > 0$ such that

$$(2.6) \qquad \langle h, DV(x)h \rangle \le -K_V$$

for $h \in \mathbb{R}^d$ s.t. $\|h\| = 1$ and $x \in \mathbb{R}^d$ s.t. $\|x\| \ge R_0$.

The conditions that make our diffusion dissipative are (2.4) and (2.6). Equation (2.4) means that the interaction is essentially not more complicated than in the one-dimensional situation and has some important implications for the geometry of the drift component $\mathbb{E}[\Phi(x - X_t^\varepsilon)]$ originating from self-interaction, namely that it points back to the origin. The same holds true for $V$ due to (2.6). In the gradient case $V = -\nabla U$, $-DV$ is the Hessian of $U$, and (2.6) means that its eigenvalues are uniformly bounded from below (w.r.t. $x$) on neighborhoods of $\infty$. Equation (2.5) is just a convenient way to combine polynomial growth and the local Lipschitz assumption in one condition. In the following two lemmas we summarize a few simple consequences of these assumptions.

LEMMA 2.2. *There exist constants $K, \eta, R_1 > 0$ such that the following hold true:*

(a) *For all $x, y \in \mathbb{R}^d$*

$$(2.7) \qquad \langle x - y, V(x) - V(y) \rangle \leq K \|x - y\|^2.$$

(b) *For $x, y \in \mathbb{R}^d$ such that $\|x - y\| \geq R_1$*

$$(2.8) \qquad \langle x - y, V(x) - V(y) \rangle \leq -\eta \|x - y\|^2.$$

(c) *For $x \in \mathbb{R}^d$ with $\|x\| \geq R_1$*

$$(2.9) \qquad \langle x, V(x) \rangle \leq -\eta \|x\|^2.$$

PROOF. Note first that, by continuity of $DV$, there exists $K > 0$ such that

$$\langle h, DV(x)h \rangle \leq K$$

holds for *all* $x$ and all $h$ of norm 1. Moreover, for $x, y \in \mathbb{R}^d$, $x \neq y$, we have

$$\frac{V(x) - V(y)}{\|x - y\|} = \int_0^1 DV(y + t(x - y)) \frac{x - y}{\|x - y\|} \, dt,$$

and therefore

$$(2.10) \qquad \left\langle \frac{x - y}{\|x - y\|}, \frac{V(x) - V(y)}{\|x - y\|} \right\rangle = \int_0^1 \langle h, DV(y + t\|x - y\|h)h \rangle \, dt,$$

where $h := \frac{x - y}{\|x - y\|}$. Since the integrand is bounded by $K$, this proves (a).

For (b), observe that the proportion of the line connecting $x$ and $y$ that lies inside $B_{R_0}(0)$ is at most $\frac{2R_0}{\|x - y\|}$. Hence

$$\left\langle \frac{x - y}{\|x - y\|}, \frac{V(x) - V(y)}{\|x - y\|} \right\rangle \leq K \frac{2R_0}{\|x - y\|} - K_V \left( 1 - \frac{2R_0}{\|x - y\|} \right),$$

which yields (b).

Part (c) is shown in a similar way. Let $x \in \mathbb{R}^d$ with $\|x\| > R_0$, and set $y := R_0 \frac{x}{\|x\|}$. Then the same argument shows the sharper bound

$$-K_V \geq \left\langle \frac{x-y}{\|x-y\|}, \frac{V(x)-V(y)}{\|x-y\|} \right\rangle = \left\langle \frac{x}{\|x\|}, \frac{V(x)-V(y)}{\|x\|-R_0} \right\rangle,$$

since the line connecting $x$ and $y$ does not intersect $B_{R_0}(0)$. Hence

$$\langle x, V(x) \rangle \leq -K_V \|x\|(\|x\|-R_0) + \|x\|\|V(y)\|,$$

which shows that (2.9) is satisfied if we set

$$R_1 = \max\left\{2R_0, 4 \sup_{\|y\|=R_0} \frac{\|V(y)\|}{K_V}\right\}$$

and $\eta = \frac{K_V}{4}$. □

LEMMA 2.3.  *For all $x, y, z \in \mathbb{R}^d$ we have:*

(a) $\|\Phi(x-y)\| \leq 2K + (K + 2^{r+1})(\|x\|^{r+1} + \|y\|^{r+1})$.
(b) $\|\Phi(x-z) - \Phi(y-z)\| \leq \|x-y\|[K + 2^r(\|x\|^r + \|y\|^r + 2\|z\|^r)]$.
(c) $\|\Phi(x-y) - \Phi(x-z)\| \leq K_1\|y-z\|(1+\|x\|^r)(1+\|y\|^r + \|z\|^r)$, *where* $K_1 = \max(K, 2^{r+1})$.

(d) *For all $x, y \in \mathbb{R}^d$ and $n \in \mathbb{N}$*

$$\langle x\|x\|^n - y\|y\|^n, \Phi(x-y) \rangle \geq 0. \tag{2.11}$$

PROOF.  By (2.5) and since $\Phi(0) = 0$ we have

$$\begin{aligned}
\|\Phi(x-y)\| &\leq \|x-y\|(K + \|x-y\|^r) \\
&\leq K(\|x\| + \|y\|) + 2^{r+1}(\|x\|^{r+1} + \|y\|^{r+1}) \\
&\leq K(2 + \|x\|^{r+1} + \|y\|^{r+1}) + 2^{r+1}(\|x\|^{r+1} + \|y\|^{r+1}) \\
&= 2K + (K + 2^{r+1})(\|x\|^{r+1} + \|y\|^{r+1}),
\end{aligned}$$

that is, (a) is proved. For (b), we use (2.5) again to see that

$$\begin{aligned}
\|\Phi(x-z) - \Phi(y-z)\| &\leq \|x-y\|(K + \|x-z\|^r + \|y-z\|^r) \\
&\leq \|x-y\|[K + 2^r(\|x\|^r + \|y\|^r + 2\|z\|^r)].
\end{aligned}$$

Property (c) follows from $\Phi(-x) = -\Phi(x)$ by further exploiting (b) as follows. We have

$$\|\Phi(x-y) - \Phi(x-z)\| \leq \|y-z\|[K + 2^{r+1}(\|x\|^r + \|y\|^r + \|z\|^r)],$$

which obviously yields (c). Finally, (d) follows from a simple calculation and (2.4). Obviously, (2.11) is equivalent to $\langle x\|x\|^n - y\|y\|^n, x-y \rangle \geq 0$. But this is an immediate consequence of the Schwarz inequality. □

Let us now return to the construction of a solution to (1.1), that is, a solution to the pair (2.1) and (2.2). The crucial property of these coupled equations is that the drift $b$ depends on (the law of) $X^\varepsilon$ and therefore also on $V$, $\varepsilon$ and the initial condition $x_0$. This means that a solution of (2.1) and (2.2) consists of a pair $(X, b)$, a continuous stochastic process $X$ and a drift term $b$, that satisfies these two equations.

Our construction of such a pair $(X, b)$ shall focus on the existence of the interaction drift $b$. It will be constructed as a fixed point in an appropriate function space such that the corresponding solution of (2.1) fulfills (2.2). Let us first derive some properties of $b$ that follow from (2.2).

LEMMA 2.4.  *Let $T > 0$, and let $(X_t)_{0 \le t \le T}$ be a stochastic process such that $\sup_{0 \le t \le T} \mathbb{E}[\|X_t\|^{r+1}] < \infty$. Then $b(t, x) = \mathbb{E}[\Phi(x - X_t)]$ has the following properties:*

(a) *$b$ is locally Lipschitz w.r.t. $x \in \mathbb{R}^d$, and the Lipschitz constant is independent of $t \in [0, T]$.*

(b) *$\langle x - y, b(t, x) - b(t, y) \rangle \ge 0$ for all $x, y \in \mathbb{R}^d$, $t \in [0, T]$.*

(c) *$b$ grows polynomially of order $r + 1$.*

PROOF.  Note first that $y \mapsto \Phi(x - y)$ grows polynomially of order $r + 1$ by Lemma 2.3(a), so that $b$ is well defined. Moreover, we have

$$\|b(t, x)\| \le \mathbb{E}[\|\Phi(x - X_t)\|] \le 2K + (K + 2^{r+1})(\|x\|^{r+1} + \mathbb{E}[\|X_t\|^{r+1}]),$$

which proves (c). For (a) observe that, by Lemma 2.3(b), we have for $z \in \mathbb{R}^d$, $x, y \in B_R(0)$,

$$\|\Phi(x - z) - \Phi(y - z)\| \le \|x - y\|[K + 2^{r+1}(R^r + \|z\|^r)].$$

Hence

$$\begin{aligned}
\|b(t, x) - b(t, y)\| &\le \mathbb{E}[\|\Phi(x - X_t) - \Phi(y - X_t)\|] \\
&\le \|x - y\|[K + 2^{r+1}(R^r + \mathbb{E}[\|X_t\|^r])]
\end{aligned}$$

for $x, y \in B_R(0)$. Since $\sup_{0 \le t \le T} \mathbb{E}[\|X_t\|^{r+1}] < \infty$, this implies (a).

In order to prove (b), fix $t \in [0, T]$, and let $\mu = \mathbb{P} \circ X_t^{-1}$. Then

$$\begin{aligned}
&\langle x - y, b(t, x) - b(t, y) \rangle \\
&\quad = \int \left\langle x - y, \frac{x - u}{\|x - u\|} \phi(\|x - u\|) - \frac{y - u}{\|y - u\|} \phi(\|y - u\|) \right\rangle \mu(du).
\end{aligned}$$

The integrand is nonnegative. Indeed, it equals

$$\|x - u\| \phi(\|x - u\|) + \|y - u\| \phi(\|y - u\|) - \left\langle y - u, \frac{x - u}{\|x - u\|} \phi(\|x - u\|) \right\rangle$$

$$-\left\langle x-u, \frac{y-u}{\|y-u\|}\phi(\|y-u\|)\right\rangle$$

$$\geq \|x-u\|\phi(\|x-u\|) + \|y-u\|\phi(\|y-u\|) - \|y-u\|\phi(\|x-u\|)$$

$$- \|x-u\|\phi(\|y-u\|)$$

$$= (\|x-u\| - \|y-u\|)(\phi(\|x-u\|) - \phi(\|y-u\|)),$$

which is nonnegative since $\phi$ is increasing, so (b) is established. $\quad\square$

In the light of the preceding lemma it is reasonable to define a space of functions that satisfy the above stated conditions, and to look for a candidate for the drift function in this space. Let $T > 0$, and for a continuous function $b : [0, T] \times \mathbb{R}^d \to \mathbb{R}^d$ define

$$(2.12) \qquad \|b\|_T := \sup_{t \in [0,T]} \sup_{x \in \mathbb{R}^d} \frac{\|b(t,x)\|}{1 + \|x\|^{2q}},$$

where $q \in \mathbb{N}$ is a fixed constant such that $2q > r$, the order of the polynomial growth of the interaction function $\Phi$. Furthermore, let

$$(2.13) \qquad \overline{\Lambda}_T := \{b : [0,T] \times \mathbb{R}^d \to \mathbb{R}^d \mid \|b\|_T < \infty, x \mapsto b(t,x)$$

$$\text{is locally Lipschitz, uniformly w.r.t. } t\}.$$

Lemma 2.4 shows that, besides being an element of $\overline{\Lambda}_T$, the drift of (2.1) must satisfy the dissipativity condition

$$(2.14) \qquad \langle x-y, b(t,x) - b(t,y)\rangle \geq 0, \qquad x, y \in \mathbb{R}^d.$$

Therefore, we define

$$(2.15) \qquad \Lambda_T := \{b \in \overline{\Lambda}_T : b \text{ satisfies } (2.14)\}.$$

It is obvious that $\| \cdot \|_T$ is indeed a norm on the vector space $\overline{\Lambda}_T$. The subset $\Lambda_T$ will be the object of interest for our construction of the interaction drift in what follows, that is, we shall construct the interaction drift as an element of $\Lambda_T$ for a proper choice of the time horizon $T$.

Once we have constructed the drift, the diffusion $X$ will simply be given as the unique strong solution of (2.1) due to the following rather classical result about strong solvability of SDEs. It ensures the existence of a unique strong solution to (2.1) for a *given* drift $b$ and is a consequence of Theorem 10.2.2 in [13], since pathwise uniqueness, nonexplosion and weak solvability imply strong solvability.

PROPOSITION 2.5. *Let* $\beta : \mathbb{R}_+ \times \mathbb{R}^d \to \mathbb{R}^d$, $(t,x) \mapsto \beta(t,x)$, *be locally Lipschitz, uniformly w.r.t.* $t \in [0,T]$ *for each* $T > 0$, *and assume that*

$$\sup_{0 \leq t \leq T} \|\beta(t,0)\| < \infty$$

*for all $T > 0$. Moreover, suppose that there exists $r_0 > 0$ such that*

$$\langle x, \beta(t, x) \rangle \leq 0 \qquad \text{for } \|x\| \geq r_0.$$

*Then the SDE*

$$dX_t = \beta(t, X_t) \, dt + \sqrt{\varepsilon} \, dW_t$$

*admits a unique strong and nonexploding solution for any random initial condition $X_0$.*

It is easily seen that the drift $\beta(t, x) = V(x) - b(t, x)$ does indeed satisfy the assumptions of Proposition 2.5 for any $b \in \Lambda_T$. This is an immediate consequence of (2.9) and (2.14).

According to Lemma 2.4, the geometric assumptions formulated in Assumption 2.1 imply that the drift term $b$ of the self-stabilizing diffusion (1.1) is an element of $\Lambda_T$, provided the moment condition stated there is satisfied. This moment condition is crucial for our construction of the drift, which motivates the following definition.

DEFINITION 2.6.   Let $T > 0$. By a *solution* of (1.1) on the time interval $[0, T]$ we mean a stochastic process $(X_t^\varepsilon)_{0 \leq t \leq T}$ that satisfies (2.1) and (2.2) on $[0, T]$ and

$$(2.16) \qquad\qquad \sup_{0 \leq t \leq T} \mathbb{E}[\|X_t\|^{2q}] < \infty.$$

A solution of (1.1) on $[0, \infty)$ is by definition a solution on $[0, T]$ for each $T > 0$.

To construct a solution of (1.1) on $[0, \infty)$, we proceed in two steps. In the first and technically most demanding step, we construct a drift on a small time interval $[0, T]$. We shall define an operator $\Gamma$ such that (2.2) translates into a fixed point property for this operator. To ensure the existence of a fixed point, one needs contraction properties of $\Gamma$ which shall turn out to depend on the time horizon $T$. This way we obtain a drift defined on $[0, T]$ such that the associated solution $X$ exists up to time $T$. In a second step, we show that this solution's moments are uniformly bounded w.r.t. time, which guarantees nonexplosion and allows us to extend $X$ to the whole time axis.

To carry out this program, we start by comparing diffusions with different drift terms.

LEMMA 2.7.   *For $b^1, b^2 \in \Lambda_T$ consider the associated diffusions*

$$dY_t = V(Y_t) \, dt - b^1(t, Y_t) \, dt + \sqrt{\varepsilon} \, dW_t$$

*and*

$$dZ_t = V(Z_t)\,dt - b^2(t, Z_t)\,dt + \sqrt{\varepsilon}\,dW_t,$$

*and assume* $Y_0 = Z_0$. *Then for* $t \le T$

$$\|Y_t - Z_t\| \le e^{KT}\|b^1 - b^2\|_T \int_0^t (1 + \|Z_s\|^{2q})\,ds.$$

PROOF. Since $Y - Z$ is governed by a (pathwise) ODE, we have

$$\|Y_t - Z_t\| = \int_0^t \left\langle \frac{Y_s - Z_s}{\|Y_s - Z_s\|}, V(Y_s) - V(Z_s) \right\rangle ds$$
$$- \int_0^t \left\langle \frac{Y_s - Z_s}{\|Y_s - Z_s\|}, b^1(s, Y_s) - b^1(s, Z_s) \right\rangle ds$$
$$+ \int_0^t \left\langle \frac{Y_s - Z_s}{\|Y_s - Z_s\|}, b^2(s, Z_s) - b^1(s, Z_s) \right\rangle ds.$$

The second integral in this decomposition is positive by definition of $\Lambda_T$, so it can be neglected. Furthermore, the first integral is bounded by $K \int_0^t \|Y_s - Z_s\|\,ds$ due to the dissipativity condition (2.7) on $V$. The last integral is bounded by

$$\int_0^t \|b^2(s, Z_s) - b^1(s, Z_s)\|\,ds \le \|b^1 - b^2\|_T \int_0^t (1 + \|Z_s\|^{2q})\,ds.$$

Combining these estimates yields

$$\|Y_t - Z_t\| \le K \int_0^t \|Y_s - Z_s\|\,ds + \|b^1 - b^2\|_T \int_0^t (1 + \|Z_s\|^{2q})\,ds.$$

Now an application of Gronwall's lemma completes the proof. □

The liberty of choice for the drift terms in Lemma 2.7 allows us to get bounds on $Y$ and its moments by making a particular one for $Z$. We consider the special case of a linear drift term $b(t, x) = \lambda x$.

LEMMA 2.8. *Let* $\lambda \ge K$, *and let* $Z$ *be the solution of*

$$dZ_t = V(Z_t)\,dt - \lambda Z_t dt + \sqrt{\varepsilon}\,dW_t.$$

*Furthermore, assume that* $\mathbb{E}(\|Z_0\|^{2m}) < \infty$ *for some* $m \in \mathbb{N}$, $m \ge 1$.
*Then for all* $t \ge 0$

$$\mathbb{E}[\|Z_t\|^{2m}] \le 2mt\|V(0)\|R_1^{2m-1}\exp\left\{\frac{\varepsilon m(d + 2m - 2)t}{R_1^2}\right\} \qquad \text{if } Z_0 = 0 \text{ a.s.,}$$

*and*

$$\mathbb{E}[\|Z_t\|^{2m}] \leq \mathbb{E}[\|Z_0\|^{2m}] \exp\left\{\frac{\varepsilon m(d+2m-2)t}{(\mathbb{E}[\|Z_0\|^{2m}])^{1/m}}\right\}$$
$$+ 2mt\|V(0)\|R_1^{2m-1}\exp\left\{\frac{\varepsilon m(d+2m-2)t}{R_1^2}\right\},$$

*otherwise.*

PROOF. By Itô's formula we have for $n \geq 2$

$$\|Z_t\|^n = \|Z_0\|^n + M_t^n + n\int_0^t \|Z_s\|^{n-2}\langle Z_s, V(Z_s)\rangle - \lambda\|Z_s\|^n\,ds$$
(2.17)
$$+ \frac{\varepsilon}{2}(dn+n^2-2n)\int_0^t \|Z_s\|^{n-2}\,ds,$$

where $M^n$ is the local martingale $M_t^n = n\sqrt{\varepsilon}\int_0^t \langle Z_s\|Z_s\|^{n-2}, dW_s\rangle$.
Since $\langle x, V(x)\rangle \leq -\eta\|x\|^2$ for $\|x\| > R_1$ according to (2.9), the first integrand of (2.17) is negative if $\|Z_s\| > R_1$. If $\|Z_s\| \leq R_1$, we use the global estimate $\langle x, V(x)\rangle \leq K\|x\|^2 + \|V(0)\|\|x\|$, which follows from (2.7). We deduce that, since $\lambda \geq K$,

$$\|Z_s\|^{n-2}\langle Z_s, V(Z_s)\rangle - \lambda\|Z_s\|^n \leq (K-\lambda)\|Z_s\|^n + \|V(0)\|\|Z_s\|^{n-1}$$
$$\leq \|V(0)\|R_1^{n-1}.$$

Thus,

$$\|Z_t\|^n \leq \|Z_0\|^n + M_t^n + n\|V(0)\|tR_1^{n-1} + \frac{\varepsilon}{2}(dn+n^2-2n)\int_0^t \|Z_s\|^{n-2}\,ds.$$

Using a localization argument and monotone convergence yields

$$\mathbb{E}[\|Z_t\|^n] \leq \mathbb{E}[\|Z_0\|^n] + n\|V(0)\|tR_1^{n-1}$$
(2.18)
$$+ \frac{\varepsilon}{2}(dn+n^2-2n)\int_0^t \mathbb{E}[\|Z_s\|^{n-2}]\,ds.$$

We claim that this implies

$$\mathbb{E}[\|Z_t\|^{2m}] \leq \sum_{j=0}^m \mathbb{E}[\|Z_0\|^{2(m-j)}]\frac{(\alpha_m t)^j}{j!}$$
(2.19)
$$+ 2m\frac{\|V(0)\|}{\alpha_m}R_1^{2m+1}\sum_{j=1}^m \frac{(\alpha_m t)^j}{R_1^{2j}j!}$$

for all $m \in \mathbb{N}$, $m \geq 1$, where $\alpha_m = \varepsilon m(d+2m-2)$. Indeed, for $m = 1$ this is evidently true by (2.18). The general case follows by induction. Assume

(2.19) holds true for $m-1$. Then by (2.18)

$$\mathbb{E}[\|Z_t\|^{2m}] \leq \mathbb{E}[\|Z_0\|^{2m}] + 2m\|V(0)\|tR_1^{2m-1}$$

$$+ \alpha_m \int_0^t \sum_{j=1}^m \mathbb{E}[\|Z_0\|^{2(m-j)}] \frac{(\alpha_{m-1}s)^{j-1}}{(j-1)!}$$

$$+ 2(m-1)\frac{\|V(0)\|}{\alpha_{m-1}} R_1^{2m-1} \sum_{j=2}^m \frac{(\alpha_{m-1}s)^{j-1}}{R_1^{2(j-1)}(j-1)!}\, ds$$

$$\leq \mathbb{E}[\|Z_0\|^{2m}] + 2m\|V(0)\|tR_1^{2m-1}$$

$$+ \sum_{j=1}^m \alpha_m \mathbb{E}[\|Z_0\|^{2(m-j)}] \frac{\alpha_{m-1}^{j-1}t^j}{j!}$$

$$+ 2m\|V(0)\|R_1^{2m-1} \sum_{j=2}^m \alpha_m \frac{\alpha_{m-1}^{j-2}t^j}{R_1^{2(j-1)}j!}$$

$$\leq 2m\|V(0)\|tR_1^{2m-1} + \sum_{j=0}^m \mathbb{E}[\|Z_0\|^{2(m-j)}] \frac{\alpha_m^j t^j}{j!}$$

$$+ 2m\|V(0)\|R_1^{2m+1} \sum_{j=2}^m \frac{\alpha_m^{j-1}t^j}{R_1^{2j}j!}$$

$$= \sum_{j=0}^m \mathbb{E}[\|Z_0\|^{2(m-j)}] \frac{\alpha_m^j t^j}{j!} + 2m\|V(0)\|R_1^{2m+1} \sum_{j=1}^m \frac{\alpha_m^{j-1}t^j}{R_1^{2j}j!},$$

and so (2.19) is established. Since $\mathbb{E}[\|Z_0\|^{2(m-j)}] \leq (\mathbb{E}[\|Z_0\|^{2m}])^{1-j/m}$ for $j \leq m$, we may exploit (2.19) further to conclude that

$$\mathbb{E}[\|Z_t\|^{2m}] \leq \mathbb{E}[\|Z_0\|^{2m}] \sum_{j=0}^m \frac{\alpha_m^j t^j}{j!(\mathbb{E}[\|Z_0\|^{2m}])^{j/m}}$$

$$+ 2mt\|V(0)\|R_1^{2m-1} \sum_{j=1}^m \frac{\alpha_m^{j-1}t^{j-1}}{R_1^{2j-2}j!}$$

$$\leq \mathbb{E}[\|Z_0\|^{2m}] \exp\left\{\frac{\alpha_m t}{(\mathbb{E}[\|Z_0\|^{2m}])^{1/m}}\right\}$$

$$+ 2mt\|V(0)\|R_1^{2m-1} \exp\left\{\frac{\alpha_m t}{R_1^2}\right\},$$

which is the announced bound if we identify the first term as zero in case $Z_0 = 0$.
$\square$

Let us define the mapping $\Gamma$ on $\Lambda_T$ that will be a contraction under suitable conditions. For $b \in \Lambda_T$, denote by $X^{(b)}$ the solution of

$$(2.20) \qquad dX_t = V(X_t)\,dt - b(t, X_t)\,dt + \sqrt{\varepsilon}\,dW_t,$$

and let $\Gamma b(t,x) := \mathbb{E}[\Phi(x - X_t^{(b)})]$. By combining the two previous lemmas, we obtain the following a priori bound on the moments of $X^{(b)}$.

LEMMA 2.9.   *If the initial datum of (2.20) satisfies $\mathbb{E}[\|X_0^{(b)}\|^{2qn}] < \infty$ for some $n \in \mathbb{N}$, then for each $T > 0$ there exists $k = k(n,T) > 0$ such that for all $b \in \Lambda_T$*

$$\sup_{0 \le t \le T} \mathbb{E}[\|X_t^{(b)}\|^n] \le k(1 + T^n e^{nKT}(\|b\|_T^n + K^n)).$$

PROOF.   Let $b^1(t,x) := b(t,x)$ and $b^2(t,x) = Kx$, and denote by $Y$, $Z$ the diffusions associated with $b^1$, $b^2$. By Lemma 2.7 we have for $t \in [0, T]$

$$\begin{aligned}
\mathbb{E}[\|Y_t\|^n] &\le 2^n(\mathbb{E}[\|Z_t\|^n] + \mathbb{E}[\|Y_t - Z_t\|^n]) \\
&\le 2^n \mathbb{E}[\|Z_t\|^n] + 2^n e^{nKT} t^{n-1} \|b^1 - b^2\|_T^n \mathbb{E}\left[\int_0^t (1 + \|Z_s\|^{2q})^n\,ds\right] \\
&\le 2^n(1 + \mathbb{E}[\|Z_t\|^{2qn}]) \\
&\qquad + 2^n e^{nKT} t^n (\|b^1\|_T + \|b^2\|_T)^n \sup_{0 \le s \le T} \mathbb{E}[(1 + \|Z_s\|^{2q})^n] \\
&\le 8^n\left(1 + \sup_{0 \le s \le T} \mathbb{E}[\|Z_s\|^{2qn}]\right)(1 + t^n e^{nKT}(\|b^1\|_T^n + \|b^2\|_T^n)).
\end{aligned}$$

Due to the assumption $\mathbb{E}[\|X_0^{(b)}\|^{2qn}] < \infty$, the constant $k(n,T) = 8^n(1 + \sup_{0 \le s \le T} \mathbb{E}[\|Z_s\|^{2qn}])$ is finite by Lemma 2.8. Furthermore, we have $\|b^2\|_T \le K$, that is, the lemma is proved.   □

Now we are in a position to establish the local Lipschitz continuity of the operator $\Gamma$. The explicit expression for the Lipschitz constant shows that $\Gamma$ will be a contraction on a sufficiently small time interval.

LEMMA 2.10.   *Let $b^1, b^2 \in \Lambda_T$, and denote by $Y, Z$ the corresponding diffusions as in Lemma 2.7. For $i \in \mathbb{N}$ let $m_i(T) = \sup_{0 \le t \le T} \mathbb{E}[\|Y_t\|^i]$ and $n_i(T) = \sup_{0 \le t \le T} \mathbb{E}[\|Z_t\|^i]$.*

*There exists a constant $k = k(m_{4q}(T), n_{4q}(T))$ such that*

$$\|\Gamma b_1 - \Gamma b_2\|_T \le k\sqrt{T}\,e^{KT}\|b^1 - b^2\|_T.$$

Proof. From Lemma 2.3(c) and the Cauchy–Schwarz inequality follows that

$$\|\Gamma b^1(t,x) - \Gamma b^2(t,x)\|$$
$$\leq \mathbb{E}[\|\Phi(x - Y_t) - \Phi(x - Z_t)\|]$$
$$\leq K_1(1 + \|x\|^r)\mathbb{E}[\|Y_t - Z_t\|(1 + \|Y_t\|^r + \|Z_t\|^r)]$$
$$\leq K_1(1 + \|x\|^r)\sqrt{\mathbb{E}[\|Y_t - Z_t\|^2]\mathbb{E}[(1 + \|Y_t\|^r + \|Z_t\|^r)^2]},$$

where $K_1 = \max(K, 2^{r+1})$. By Lemma 2.7, since $(1 + x)^2 \leq 2(1 + x^2)$, we have

$$\mathbb{E}[\|Y_t - Z_t\|^2] \leq e^{2KT}\|b^1 - b^2\|_T^2 \mathbb{E}\left[\left(\int_0^T (1 + \|Z_s\|^{2q})\,ds\right)^2\right]$$
$$\leq e^{2KT}\|b^1 - b^2\|_T^2 \int_0^T \mathbb{E}[(1 + \|Z_s\|^{2q})^2]\,ds$$
$$\leq 2Te^{2KT}\|b^1 - b^2\|_T^2\left(1 + \sup_{0\leq s\leq T}\mathbb{E}[\|Z_s\|^{4q}]\right).$$

Moreover, using the inequality $(a + b)^2 \leq 2(a^2 + b^2)$, we deduce that

$$\mathbb{E}[(1 + \|Y_t\|^r + \|Z_t\|^r)^2] \leq 2(1 + 2\mathbb{E}[\|Y_t\|^{2r} + \|Z_t\|^{2r}])$$
$$\leq 10(1 + \mathbb{E}[\|Y_t\|^{4q} + \|Z_t\|^{4q}]),$$

where we exploited that $2q > r$ implies $\mathbb{E}[\|Y_t\|^{2r}] \leq 1 + \mathbb{E}[\|Y_t\|^{4q}]$, and likewise for the moment of $Z_t$. By combining all these estimates, we find that

$$\frac{\|\Gamma b^1(t,x) - \Gamma b^2(t,x)\|}{1 + \|x\|^{2q}}$$
$$\leq 2K_1\sqrt{5T}e^{KT}\|b^1 - b^2\|_T \frac{1 + \|x\|^r}{1 + \|x\|^{2q}}$$
$$\times \left(1 + \sup_{0\leq s\leq T}\mathbb{E}[\|Z_s\|^{4q}]\right)^{1/2}(1 + \mathbb{E}[\|Y_t\|^{4q} + \|Z_t\|^{4q}])^{1/2}.$$

Hence, if we set $k := 4K_1\sqrt{5}\{(1 + n_{4q}(T))(1 + m_{4q}(T) + n_{4q}(T))\}^{1/2}$, we may conclude that

$$\|\Gamma b^1 - \Gamma b^2\|_T \leq k\sqrt{T}e^{KT}\|b^1 - b^2\|_T,$$

that is, $k$ is the desired constant. $\square$

The next proposition shows that the restriction of $\Gamma$ to a suitable subset of the function space $\Lambda_T$ is a contractive mapping, which allows us to construct a solution on a small time interval.

PROPOSITION 2.11.    *For $\nu > 0$ let $\Lambda_T^\nu = \{b \in \Lambda_T : \|b\|_T \leq \nu\}$. Assume that the initial condition $X_0$ satisfies $\mathbb{E}[\|X_0\|^{2qn}] < \infty$ for some $n \geq 4q$. There exists $\nu_0 > 0$ such that for any $\nu \geq \nu_0$ there exists $T = T(\nu) > 0$ such that the following hold true:*

(a) *$\Gamma(\Lambda_T^\nu) \subset \Lambda_T^\nu$, and the Lipschitz constant of $\Gamma|\Lambda_T^\nu$ is less than $\frac{1}{2}$.*
(b) *There exists a strong solution to (2.1), (2.2) on $[0, T]$ which satisfies*

$$\sup_{0 \leq t \leq T} \mathbb{E}[\|X_t^{(b)}\|^n] \leq k(1 + T^n e^{nKT}(\nu^n + K^n)),$$

*where $k = k(n, T)$ is the constant introduced in Lemma 2.9.*

PROOF.    Let $b \in \Lambda_T$, and let $X = X^{(b)}$ and $m_i(T) = \sup_{0 \leq t \leq T} \mathbb{E}[\|X_t\|^i]$ for $i \in \mathbb{N}$. By Lemma 2.9 the condition $\mathbb{E}[\|X_0\|^{2qn}] < \infty$ implies $m_i(T) < \infty$ for $T > 0$ and $i \leq n$. Moreover, Lemma 2.3 shows that

$$\|\Gamma b(t, x)\| \leq 2K + (K + 2^{r+1})(\|x\|^{r+1} + \mathbb{E}[\|X\|^{r+1}])$$
$$\leq \tilde{K}(1 + \|x\|^{r+1})(1 + \mathbb{E}[\|X_t\|^{r+1}]),$$

where $\tilde{K} = 2K + 2^{r+1}$. Consequently, by definition of $\|\cdot\|_T$,

$$(2.21) \qquad \|\Gamma b\|_T \leq 2\tilde{K}(1 + m_{r+1}(T)), \qquad t \leq T.$$

By Lemma 2.9 there exists $k = k(r + 1, T) > 0$ such that

$$(2.22) \qquad m_{r+1}(T) \leq k(1 + T^{r+1} e^{(r+1)KT}(\|b\|_T^{r+1} + K^{r+1})).$$

This inequality, together with (2.21), is the key for finding a suitable subset of $\Lambda_T$ on which $\Gamma$ is contractive. The r.h.s. of (2.22) converges to $k$ as $T \to 0$, and this convergence is uniform w.r.t. $b \in \Lambda_T^\nu$ for each $\nu > 0$. The dependence of the limiting constant $k$ on $T$ imposes no problem here; just fix $k = k(r + 1, T_0) > 0$ for some $T_0$ and use the fact that (2.22) is valid for all $T \leq T_0$, as the proof of Lemma 2.9 shows.

Thus, we may fix $\nu_0 > 2\tilde{K}(1 + k)$ and deduce that for any $\nu > \nu_0$ we can find $T_0 = T_0(\nu)$ such that $\|b\|_T \leq \nu$ implies $\|\Gamma b\|_T \leq \nu$ for $T \leq T_0$. Moreover, by Lemma 2.4, $\Gamma b$ satisfies all the conditions as required for it to belong to $\Lambda_T$, that is, $\Gamma$ maps $\Lambda_T^\nu$ into itself for all $T \leq T_0$. Additionally, the assumption $n \geq 4q$ implies that $m_{4q}(T)$ is uniformly bounded for all $b$ in $\Lambda_T^\nu$, and Lemma 2.10 shows that, by eventually decreasing $T_0$, we can achieve that $\Gamma$ is a contraction on $\Lambda_T^\nu$ with Lipschitz constant less than $\frac{1}{2}$, that is, (a) is established.

In order to prove (b), the existence of a strong solution on the time interval $[0, T]$ for some $T \leq T_0$, we iterate the drift through $\Gamma$. Let $b_0 \in \Lambda_T^\nu$, and define

$$b_{i+1} := \Gamma b_i \qquad \text{for } i \in \mathbb{N}_0.$$

The contraction property of $\Gamma$ yields $\|b_{i+1} - b_i\|_T \leq 2^{-i}\|b_1 - b_0\|_T$ for all $i$, and therefore

$$\sum_{n=0}^{\infty} \|b_{i+1} - b_i\|_T < \infty,$$

which entails that $(b_i)$ is a Cauchy sequence w.r.t. $\|\cdot\|_T$. By definition of $\|\cdot\|_T$, $(b_i)$ converges pointwise to a continuous function $b = b(t, x)$ with $\|b\|_T < \infty$. It remains to verify that the limit is again an element of $\Lambda_T$. In order to see that it is locally Lipschitz, let $X^{(i)} := X^{(b_i)}$. As in the proof of Lemma 2.4, we have for $x, y \in B_R(0)$

$$\|\Gamma b_i(t, x) - \Gamma b_i(t, y)\| \leq \mathbb{E}[\|\Phi(x - X_t^{(i)}) - \Phi(y - X_t^{(i)})\|]$$
$$\leq \|x - y\|[K + 2^{r+1}(R^r + \mathbb{E}[\|X_t^{(i)}\|^r])].$$

Since $\|b_i\|_T \leq \nu$ for all $i$, (2.22) yields

$$\sup_{i \in \mathbb{N}} \sup_{0 \leq t \leq T} \mathbb{E}[\|X_t^{(i)}\|^r] \leq k(1 + T^{r+1}e^{(r+1)KT}(\nu^{r+1} + K^{r+1})).$$

Therefore, we may send $i \to \infty$ to conclude that $b$ is locally Lipschitz. $b$ being the pointwise limit of the $b_i$, it inherits the polynomial growth property and the dissipativity condition as stated in Lemma 2.4(b) and (c). (Notice that we may not invoke Lemma 2.4 at this stage.)

It remains to show that the diffusion $X = X^{(b)}$ associated to $b$ has the desired properties. Note first that the existence of $X$ is guaranteed by the classical result of Proposition 2.5. Since $\Gamma b = b$, which means that

$$b(t, x) = \Gamma b(t, x) = \mathbb{E}[\Phi(x - X_t^{(b)})]$$

for $t \in [0, T]$ and $x \in \mathbb{R}^d$, $X$ is the diffusion with interaction drift $b$. The boundedness of its moments is again a consequence of Lemma 2.9. $\quad\square$

Let us recall the essentials of the construction carried out so far. We have shown the existence of a solution to (1.1) on a small time interval $[0, T]$. For the moments of order $n$ to be finite, one needs integrability of order $2qn$ for the initial condition. Moreover, the parameter $n$ needs to be larger than or equal to $4q$ in order for the fixed point argument of proposition 2.11 to work. Observe that the condition $n \geq 4q$ appears first in this proposition, since this is the first time the process is coupled to its own drift, while in all previous statements the finiteness of moments is guaranteed by the comparison against the diffusion $Z$, which is governed by a linear drift term. In order to find a solution that exists for all times, we need to carefully extend the constructed pair $(X, b)$ beyond the time horizon $T$. Although nonexplosion and finiteness of moments would be guaranteed for all $T$ by Proposition 2.5 and Lemma 2.9, we have to take care of the fact that the drift itself is defined only on the time interval $[0, T]$. With sufficiently strong integrability assumptions for $X_0$ one could perform the same construction on the time intervals $[T, 2T]$, $[2T, 3T]$ and so on, but one loses an integrability order $2q$ in each time step of length $T$.

For that reason we need better control of the moments of $X$ over the whole time axis, which is achieved by the following a posteriori estimate.

PROPOSITION 2.12.   *Let* $m \in \mathbb{N}$, $m \geq 4q^2$, *such that* $\mathbb{E}[\|X_0\|^{2m}] < \infty$. *For each* $n \in \{1, \ldots, m\}$ *there exists a constant* $\alpha = \alpha(n) > 0$ *such that the following holds true for all* $T > 0$: *if* $X$ *solves (1.1) on* $[0, T]$, *then*

$$\sup_{0 \leq t \leq T} \mathbb{E}[\|X_t\|^{2n}] \leq \alpha(n).$$

PROOF.   Fix $T > 0$, and assume that $X$ solves (1.1) on $[0, T]$ (in the sense of Definition 2.6). Then $b(t, x) := \mathbb{E}[\Phi(x - X_t)]$ belongs to $\Lambda_T$ by Lemma 2.4, and $m \geq 4q^2$ implies $\sup_{0 \leq t \leq T} \mathbb{E}[\|X_t\|^{4q}] < \infty$ by Lemma 2.9. Let $f_n(t) = \mathbb{E}[\|X_t\|^{2n}]$. We proceed in several steps.

*Step* 1: Boundedness in $L^2$. We know already (by Lemma 2.9) that

$$\sup_{0 \leq t \leq T} f_1(t) < \infty.$$

The only point is to show that the bound may be chosen independent of $T$. By Itô's formula we have

$$f_1(t) = \mathbb{E}[\|X_0\|^2] + \varepsilon t d + 2 \int_0^t \mathbb{E}[\langle X_s, V(X_s) \rangle] \, ds - 2 \int_0^t \mathbb{E}[\langle X_s, b(s, X_s) \rangle] \, ds.$$

Let us estimate the last term that contains the interaction drift $b$. Note first that $2q > r$ implies $r + 2 \leq 4q$, so $\sup_{0 \leq t \leq T} \mathbb{E}[\|X_t\|^{r+2}] < \infty$ as pointed out at the beginning of the proof, and the Cauchy–Schwarz inequality yields that $\mathbb{E}[\langle X_t, b(t, X_t) \rangle]$ is finite for $t \in [0, T]$ since $b$ grows polynomially of order $r + 1$. By definition of $b$, we may take an independent copy $\tilde{X}$ of $X$, to write

$$\begin{aligned} 2\mathbb{E}[\langle X_s, b(s, X_s) \rangle] &= 2\mathbb{E}[\langle X_s, \Phi(X_s - \tilde{X}_s) \rangle] \\ &= \mathbb{E}[\langle X_s, \Phi(X_s - \tilde{X}_s) \rangle] - \mathbb{E}[\langle \tilde{X}_s, \Phi(X_s - \tilde{X}_s) \rangle] \\ &= \mathbb{E}[\langle X_s - \tilde{X}_s, \Phi(X_s - \tilde{X}_s) \rangle] \geq 0 \end{aligned}$$

where the last inequality is due to (2.4). In order to estimate the other integral, let $R \geq R_1$. Using (2.9) and the local Lipschitz property of $V$, we see that

$$\begin{aligned} \mathbb{E}[\langle X_s, V(X_s) \rangle] &\leq -\eta \mathbb{E}[\|X_s\|^2 \mathbf{1}_{\{\|X_s\| > R\}}] \\ &\quad + \mathbb{E}[(K\|X_s\|^2 + \|V(0)\| \|X_s\|) \mathbf{1}_{\{\|X_s\| \leq R\}}] \end{aligned}$$

$$\leq -\eta \mathbb{E}[\|X_s\|^2] + (\eta + K)R^2 + \|V(0)\|R$$

$$= -\eta f_1(s) + R(\|V(0)\| + R(\eta + K)).$$

Obviously, $f_1$ is differentiable, and summing up these bounds yields

$$f_1'(t) \leq \varepsilon d - 2\eta f_1(t) + 2R(\|V(0)\| + R(\eta + K)).$$

Thus, there exists $\gamma > 0$ such that $\{t \in [0,T] : f_1(t) \geq \gamma\} \subset \{t \in [0,T] : f_1'(t) \leq 0\}$, which implies $f_1(t) \leq f_1(0) \vee \gamma$ for all $t \in [0,T]$. This is the claimed bound, since $\gamma$ is independent of $T$.

*Step* 2: Moment bound for the convolution. Let $\tilde{X}$ be an independent copy of $X$, that is, a solution of (1.1) driven by a Brownian motion that is independent of $W$. In this step we shall prove that $\mathbb{E}[\|X_t - \tilde{X}_t\|^{2n}]$ is uniformly bounded w.r.t. time.

Let $R \geq R_1$, and let $\tau = \inf\{t \geq 0 : \|X_t - \tilde{X}_t\| \geq R\}$, $g_n(t) = \mathbb{E}[\|X_t - \tilde{X}_t\|^{2n}\mathbf{1}_{\{t < \tau\}}]$ and $w_n(t) = \mathbb{E}[\|X_{t \wedge \tau} - \tilde{X}_{\tau \wedge t}\|^{2n}]$. Then $w_n(t) = g_n(t) + R^{2n}\mathbb{P}(t \geq \tau)$. Furthermore, using the SDE (1.1) for both $X$ and $\tilde{X}$, applying Itô's formula to the difference and taking expectations, we obtain for $n \geq 1$

$$w_n(t) = \mathbb{E}[\|X_0 - \tilde{X}_0\|^{2n}] + \varepsilon n(d + 2n - 2)\mathbb{E}\left[\int_0^{t \wedge \tau} \|X_s - \tilde{X}_s\|^{2n-2} ds\right]$$

$$+ 2n\mathbb{E}\left[\int_0^{t \wedge \tau} \|X_s - \tilde{X}_s\|^{2n-2}\langle X_s - \tilde{X}_s, V(X_s) - V(\tilde{X}_s)\rangle ds\right]$$

$$- 2n\mathbb{E}\left[\int_0^{t \wedge \tau} \|X_s - \tilde{X}_s\|^{2n-2}\langle X_s - \tilde{X}_s, b(s, X_s) - b(s, \tilde{X}_s)\rangle ds\right].$$

The last term is negative by Lemma 2.4, which yields together with (2.7), (2.8) and Hölder's inequality

$$w_n'(t) \leq \varepsilon n(d + 2n - 2)\mathbb{E}[\|X_t - \tilde{X}_t\|^{2n-2}\mathbf{1}_{\{t < \tau\}}]$$

$$+ 2n\mathbb{E}[\|X_t - \tilde{X}_t\|^{2n-2}\langle X_t - \tilde{X}_t, V(X_t) - V(\tilde{X}_t)\rangle\mathbf{1}_{\{t < \tau\}}]$$

$$\leq \varepsilon n(d + 2n - 2)g_{n-1}(t)$$

$$+ 2n(K + \eta)\mathbb{E}[\|X_t - \tilde{X}_t\|^{2n}\mathbf{1}_{\{\|X_t - \tilde{X}_t\| \leq R_1; \tau > t\}}]$$

$$- 2n\eta\mathbb{E}[\|X_t - \tilde{X}_t\|^{2n}\mathbf{1}_{\{t < \tau\}}]$$

$$\leq \varepsilon n(d + 2n - 2)g_n(t)^{1-1/n} + 2n(K + \eta)R_1^{2n} - 2n\eta g_n(t).$$

As in the first step, there exists some constant $\delta > 0$ such that $\{t \in [0,T] : g_n(t) > \delta\} \subset \{t \in [0,T] : w_n'(t) < 0\}$. Since $w_n - g_n$ is nondecreasing this implies $g_n(t) \leq g_n(0) \vee \delta$ for all $t \in [0,T]$. Moreover, $\delta$ depends only on the constants appearing in the last inequality and is independent of the localization parameter $R$. Hence, by monotone convergence, we have

$$\mathbb{E}[\|X_t - \tilde{X}_t\|^{2n}] \leq \mathbb{E}[\|X_0 - \tilde{X}_0\|^{2n}] \vee \delta, \qquad t \in [0,T].$$

*Step* 3: Bound for the centered moments of $X$. In this step we shall prove that the moments of $Y_t := X_t - \mathbb{E}[X_t]$ are uniformly bounded. We proceed by induction. The second moments of $X$ are uniformly bounded by the first step; so are those of $Y$. Assume the moments of order $2n$ are uniformly bounded by $\gamma_n > 0$. If $n + 1 \leq m$, we may invoke step 2, to find $\delta_{n+1} > 0$ such that $\mathbb{E}[\|X_t - \tilde{X}_t\|^{2n+2}] \leq \delta_{n+1}$ for $t \in [0, T]$. Now we make the following observation. If $\xi$, $\tilde{\xi}$ are independent, real-valued copies of each other with $\mathbb{E}[\xi] = 0$, then

$$\mathbb{E}[(\xi - \tilde{\xi})^{2n+2}] = 2\mathbb{E}[\xi^{2n+2}] + \sum_{k=2}^{2n} \binom{2n+2}{k} (-1)^k \mathbb{E}[\xi^k] \mathbb{E}[\xi^{2n+2-k}],$$

and therefore

$$2\mathbb{E}[\xi^{2n+2}] \leq \mathbb{E}[(\xi - \tilde{\xi})^{2n+2}] + \sum_{k=2}^{2n} \binom{2n+2}{k} |\mathbb{E}[\xi^k] \mathbb{E}[\xi^{2n+2-k}]|$$

$$\leq \mathbb{E}[(\xi - \tilde{\xi})^{2n+2}] + 2^{2n+2}(1 + \mathbb{E}[\xi^{2n}])^2.$$

Let us apply this to the components of $Y$, and denote them by $Y^1, \dots, Y^d$. We obtain for $t \in [0, T]$

$$2\mathbb{E}[\|Y_t\|^{2n+2}] \leq 2d^{n+1}\mathbb{E}\left[\sum_{j=1}^d (Y_t^j)^{2n+2}\right]$$

$$\leq d^{n+1}\sum_{j=1}^d \mathbb{E}[(X_t^j - \tilde{X}_t^j)^{2n+2}] + 2^{2n+2}(1 + \mathbb{E}[(Y_t^j)^{2n}])^2$$

$$\leq d^{n+2}(\mathbb{E}[\|X_t - \tilde{X}_t\|^{2n+2}] + 2^{2n+2}(1 + \mathbb{E}[\|Y_t\|^{2n}])^2)$$

$$\leq d^{n+2}(\delta_{n+1} + 2^{2n+2}(1 + \gamma_n)^2),$$

which is a uniform bound for the order $2(n + 1)$.

*Step* 4: Bound for the moments of $X$. In the fourth and final step, we prove the announced uniform bound for the moments of $X$. It follows immediately from the inequality

$$\mathbb{E}[\|X_t\|^{2n}] \leq 2^{2n}(\mathbb{E}[\|X_t - \mathbb{E}[X_t]\|^{2n}] + \|\mathbb{E}[X_t]\|^{2n}).$$

The last term satisfies $\|\mathbb{E}[X_t]\|^{2n} \leq f_1(t)^n$, which is uniformly bounded according to step 1, and the centered moments of order $2n$ are uniformly bounded by step 3 whenever $n \leq m$.  $\square$

The results concerning the existence of $X^\varepsilon$ are summarized in the following theorem.

THEOREM 2.13. *Let $q := [\frac{r}{2} + 1]$, and let $X_0$ be a random initial condition such that $\mathbb{E}[\|X_0\|^{8q^2}] < \infty$. Then there exists a drift term $b(t, x) = b^{\varepsilon, X_0}(t, x)$ such that (2.1) admits a unique strong solution $X^{\varepsilon}$ that satisfies (2.2), and $X^{\varepsilon}$ is the unique strong solution of (1.1). Moreover, we have for all $n \in \mathbb{N}$*

$$(2.23) \qquad \sup_{t \geq 0} \mathbb{E}[\|X_t^{\varepsilon}\|^{2n}] < \infty$$

*whenever $\mathbb{E}[\|X_0^{\varepsilon}\|^{2n}] < \infty$. In particular, if $X_0$ is deterministic, then $X^{\varepsilon}$ is bounded in $L^p(\mathbb{P} \otimes \boldsymbol{\lambda}_{[0,T]})$ for all $p \geq 1$. $\boldsymbol{\lambda}$ is used as a symbol for Lebesgue measure throughout.*

PROOF. In a first step, we prove uniqueness on a small time interval. Let $\tilde{K} = 2K + 2^{r+1}$, and choose $\alpha(q) > 0$ according to Proposition 2.12. By Proposition 2.11 there exist $\nu \geq 2\tilde{K}(2 + \alpha(q))$, $T = T(\nu) > 0$ and $b \in \Lambda_T^{\nu}$ such that $\Gamma b = b$, that is, $X = X^{(b)}$ is a strong solution of (1.1) on $[0, T]$. Assume $Y$ is another solution of (1.1) on $[0, T]$ starting at $X_0$ such that $m_{2q}(T) := \sup_{0 \leq t \leq T} \mathbb{E}[\|Y_t\|^{2q}] < \infty$, and let $c(t, x) = \mathbb{E}[\Phi(x - Y_t)]$. Then $c \in \Lambda_T$ by Lemma 2.4, and $\Gamma c = c$. Moreover, it follows from (2.21) and Proposition 2.12 that

$$\|c\|_T \leq 2\tilde{K}(2 + m_{2q}(T)) \leq 2\tilde{K}(2 + \alpha(q)) \leq \nu,$$

that is, $c \in \Lambda_T^{\nu}$. Hence $c$ is the unique fixed point of $\Gamma|\Lambda_T^{\nu}$. Thus $c = b$, and Proposition 2.5 yields $X = Y$.

In the second step, we show the existence of a unique solution on $[0, \infty)$. Let

$$U := \sup \Big\{ T > 0 : (1.1) \text{ admits a unique strong solution } X \text{ on } [0, T],$$

$$\sup_{0 \leq t \leq T} \mathbb{E}[\|X_t\|^{2q}] < \infty \Big\}.$$

By the first step we know that $U > 0$. Assume $U < \infty$. As in the first step, choose $\alpha(4q^2) > 0$ according to Proposition 2.12, and then fix $\tilde{\nu} \geq 2\tilde{K}(2 + \alpha(4q^2))$ and $\tilde{T} = \tilde{T}(\tilde{\nu}) > 0$ that satisfy Proposition 2.11. Let $0 < \delta < \min(U, \tilde{T}/2)$, and fix $T \in ]U - \delta, U[$. There exists a unique strong solution $X$ on $[0, T]$, and $\mathbb{E}[\|X_T\|^{8q^2}] < \infty$ by Proposition 2.12. Now consider (1.1) on $[T, \infty)$ with initial datum $X_T$. As in the first step, we may find a unique strong solution on $[T, T + \tilde{T}]$. But this is a contradiction since $T + \tilde{T} > U$. Consequently, $U = \infty$, and (2.23) holds by Proposition 2.12. $\square$

**3. Large deviations.** Let us now turn to the large deviations behavior of the diffusion $X^\varepsilon$ given by the SDE (1.1), that is,

$$(3.1) \qquad dX_t^\varepsilon = V(X_t^\varepsilon)\,dt - \int_{\mathbb{R}^d} \Phi(X_t^\varepsilon - x)\,du_t^\varepsilon(x)\,dt + \sqrt{\varepsilon}\,dW_t, t \geq 0,$$

$$X_0 = x_0 \in \mathbb{R}^d.$$

The heuristics underlying large deviations theory is to identify a deterministic path around which the diffusion is concentrated with overwhelming probability, so that the stochastic motion can be seen as a small random perturbation of this deterministic path. This means in particular that the law $u_t^\varepsilon$ of $X_t^\varepsilon$ is close to some Dirac mass if $\varepsilon$ is small. We therefore proceed in two steps toward the aim of proving a large deviations principle for $X^\varepsilon$. In a first step we "guess" the deterministic limit around which $X^\varepsilon$ is concentrated for small $\varepsilon$, and replace $u_t^\varepsilon$ by its suspected limit, that is, we approximate the law of $X^\varepsilon$. This way we circumvent the difficulty of the dependence on the law of $X^\varepsilon$—the self-interaction term—and obtain a diffusion which is defined by means of a classical SDE. We then prove in the second step that this diffusion is exponentially equivalent to $X^\varepsilon$, that is, it has the same large deviations behavior. This involves pathwise comparisons.

3.1. *Small noise asymptotics of the interaction drift.* The limiting behavior of the diffusion $X^\varepsilon$ can be guessed in the following way. As explained, the laws $u_t^\varepsilon$ should tend to a Dirac measure in the small noise limit, and since $\Phi(0) = 0$ the interaction term will vanish in the limiting equation. Therefore, the diffusion $X^\varepsilon$ is a small random perturbation of the deterministic motion $\psi$, given as the solution of the deterministic equation

$$(3.2) \qquad\qquad \dot{\psi}_t = V(\psi_t), \qquad \psi_0 = x_0,$$

and the large deviations principle will describe the asymptotic deviation of $X^\varepsilon$ from this path. Much like in the case of gradient type systems, the dissipativity condition (2.9) guarantees nonexplosion of $\psi$. Indeed, since $\frac{d}{dt}\|\psi_t\|^2 = 2\langle \psi_t, \dot{\psi}_t \rangle = 2\langle \psi_t, V(\psi_t) \rangle$, the derivative of $\|\psi_t\|^2$ is negative for large values of $\|\psi_t\|$ by (2.9), so $\psi$ is bounded. In the sequel we shall write $\psi_t(x_0)$ if we want to stress the dependence on the initial condition.

We have to control the diffusion's deviation from this deterministic limit on a finite time interval. An a priori estimate is provided by the following lemma, which gives an $L^2$-bound for this deviation. For notational convenience, we suppress the $\varepsilon$-dependence of the diffusion in the sequel, but keep in mind that all processes depend on $\varepsilon$.

LEMMA 3.1. *Let $Z_t := X_t - \psi_t(x_0)$. Then*

$$\mathbb{E}\|Z_t\|^2 \leq \varepsilon t\,de^{2Kt},$$

*where $K$ is the constant introduced in Lemma 2.2. In particular, $Z \to 0$ as $\varepsilon \to 0$ in $L^p(\mathbb{P} \otimes \boldsymbol{\lambda}_{[0,T]})$ for all $p \geq 1$ and $T > 0$. This convergence is locally uniform w.r.t. the initial condition $x_0$.*

PROOF. By Itô's formula we have

$$\|Z_t\|^2 = 2\sqrt{\varepsilon} \int_0^t \langle Z_s, dW_s \rangle - 2 \int_0^t \langle Z_s, b^{\varepsilon, x_0}(s, Z_s + \psi_s(x_0)) \rangle \, ds$$
$$+ 2 \int_0^t \langle Z_s, V(Z_s + \psi_s(x_0)) - V(\psi_s(x_0)) \rangle \, ds + \varepsilon t d.$$

By Theorem 2.13 $X$ and thus $Z$ is integrable of all orders. In particular, $Z$ is square-integrable, so the stochastic integral in this equation is a martingale. Now consider the second term containing the interaction drift $b^{\varepsilon, x_0}$. Let $\nu_s = \mathbb{P} \circ Z_s^{-1}$ denote the law of $Z_s$. Since $Z$ has finite moments of all orders, Lemma 2.3 implies $\int \int \|\langle z, \Phi(z-y) \rangle \| \nu_s(dy) \nu_s(dz) < \infty$. Thus, by Assumption 2.1(ii) about the interaction function $\Phi$ and Fubini's theorem,

$$2\mathbb{E} \langle Z_s, b^{\varepsilon, x_0}(s, Z_s + \psi_s(x_0)) \rangle = 2 \int \langle z, \mathbb{E}[\Phi(z + \psi_s(x_0) - X_s)] \rangle \nu_s(dz)$$
$$= 2 \int \int \langle z, \Phi(z-y) \rangle \nu_s(dy) \nu_s(dz)$$
$$= \int \int \langle z-y, \Phi(z-y) \rangle \nu_s(dy) \nu_s(dz) \geq 0.$$

Hence by the growth condition (2.7) for $V$

$$\mathbb{E}\|Z_t\|^2 \leq 2 \int_0^t \mathbb{E} \langle Z_s, V(Z_s + \psi_s(x_0)) - V(\psi_s(x_0)) \rangle \, ds + \varepsilon t d$$
$$\leq 2K \int_0^t \mathbb{E}\|Z_s\|^2 \, ds + \varepsilon t d,$$

and Gronwall's lemma yields

$$\mathbb{E}\|Z_t\|^2 \leq \varepsilon t \, d \, e^{2Kt}.$$

This is the claimed bound. For the $L^p$-convergence observe that this bound is independent of the initial condition $x_0$. Moreover, the argument of Proposition 2.12 shows that $\sup\{\mathbb{E}(\|X_t\|^p) : 0 \leq t \leq T, x_0 \in L, 0 < \varepsilon < \varepsilon_0\} < \infty$ holds for compact sets $L$ and $\varepsilon_0 > 0$. This implies that $Z$ is bounded in $L^p(\mathbb{P} \otimes \boldsymbol{\lambda}_{[0,T]})$ as $\varepsilon \to 0$, uniformly w.r.t. $x_0 \in L$. Now the $L^p$-convergence follows from the Vitali convergence theorem. $\square$

COROLLARY 3.2. *For any $T > 0$ we have*

$$\lim_{\varepsilon \to 0} b^{\varepsilon, x_0}(t, x) = \Phi(x - \psi_t(x_0)),$$

*uniformly w.r.t. $t \in [0, T]$ and w.r.t. $x$ and $x_0$ on compact subsets of $\mathbb{R}^d$.*

PROOF.   The growth condition on $\Phi$ and the Cauchy–Schwarz inequality yield

$$\|b^\varepsilon(t,x) - \Phi(x - \psi_t(x_0))\|^2$$
$$\leq \mathbb{E}[\|X_t - \psi_t(x_0)\|(K + \|x - X_t\|^r + \|x - \psi_t(x_0)\|^r)]^2$$
$$\leq \mathbb{E}[\|X_t - \psi_t(x_0)\|^2]\mathbb{E}[(K + \|x - X_t\|^r + \|x - \psi_t(x_0)\|^r)^2].$$

The first expectation on the r.h.s. of this inequality tends to zero by Lemma 3.1. Since $X$ is bounded in $L^{2r}(\mathbb{P})$, uniformly w.r.t. $x_0$ on compact sets, the claimed convergence follows.   $\square$

In a next step we replace the diffusion's law in (3.1) by its limit, the Dirac measure in $\psi_t(x_0)$. Before doing so, let us introduce a slight generalization of $X$.

Theorem 2.13 implies that $X$ is a time-inhomogeneous Markov process. The diffusion $X$, starting at time $s \geq 0$, is given as the unique solution of the stochastic integral equation

$$X_t = X_s + \int_s^t [V(X_u) - b^{\varepsilon,x_0}(u,X_u)]\,du + \sqrt{\varepsilon}(W_t - W_s), \qquad t \geq s.$$

By shifting the starting time back to the origin, this equation translates into

$$X_{t+s} = X_s + \int_0^t [V(X_{u+s}) - b^{\varepsilon,x_0}(u+s,X_{u+s})]\,du + \sqrt{\varepsilon}W_t^s, \qquad t \geq 0,$$

where $W^s$ is the Brownian motion given by $W_t^s = W_{t+s} - W_s$, which is independent of $X_s$. Since we are mainly interested in the law of $X$, we may replace $W^s$ by $W$.

For an initial condition $\xi_0 \in \mathbb{R}^d$ and $s \geq 0$, we denote by $\xi^{s,\xi_0}$ the unique solution of the equation

$$(3.3) \qquad \xi_t = \xi_0 + \int_0^t V(\xi_u) - b^{\varepsilon,x_0}(u+s,\xi_u)\,du + \sqrt{\varepsilon}W_t, \qquad t \geq 0.$$

Note that $\xi^{0,x_0} = X$, and that $\xi^{s,\xi_0}$ has the same law as $X_{t+s}$, given that $X_s = \xi_0$. The interpretation of $b^{\varepsilon,x_0}$ as an interaction drift is lost in this equation, since $b^{\varepsilon,x_0}$ does not depend on $\xi^{s,\xi_0}$.

Now recall that $b^{\varepsilon,x_0}(t,x) = \mathbb{E}\{\Phi(x - X_t^\varepsilon)\}$, which tends to $\Phi(x - \psi_t(x_0))$ by Corollary 3.2. This motivates the definition of the following analogue of $\xi^{s,\xi_0}$, in which $u_t^\varepsilon$ is replaced by the Dirac measure in $\psi_t(x_0)$. We denote by $Y^{s,y}$ the solution of the equation

$$(3.4) \quad Y_t = y + \int_0^t V(Y_u) - \Phi(Y_u - \psi_{t+s}(x_0))\,du + \sqrt{\varepsilon}W_t, \qquad t \geq 0.$$

This equation is an SDE in the classical sense, and it admits a unique strong solution by Proposition 2.5. Furthermore, it is known that $Y^{s,y}$ satisfies a large deviations principle in the space

$$C_{0T} = \{f : [0, T] \to \mathbb{R}^d | f \text{ is continuous}\},$$

equipped with the topology of uniform convergence. This LDP describes the deviations of $Y^{s,y}$ from the deterministic system $\dot{\varphi}_t = V(\varphi_t) - \Phi(\varphi_t - \psi_{t+s}(x_0))$ with $\varphi_0 = y$. Observe that $\varphi$ coincides with $\psi(x_0)$ in case $y = x_0$, and that nonexplosion of $\varphi$ is ensured by the dissipativity properties of $V$ and $\Phi$ as follows. By (2.4) we have

$$\begin{aligned}
\frac{d}{dt}\|\varphi_t - \psi_{t+s}\|^2 &= 2\langle \varphi_t - \psi_{t+s}, \dot{\varphi}_t - \dot{\psi}_{t+s}\rangle \\
&= 2\langle \varphi_t - \psi_{t+s}, V(\varphi_t) - \Phi(\varphi_t - \psi_{t+s}) - V(\psi_{t+s})\rangle \\
&\leq 2\langle \varphi_t - \psi_{t+s}, V(\varphi_t) - V(\psi_{t+s})\rangle.
\end{aligned}$$

(3.5)

Since the last expression is negative for large values of $\|\varphi_t - \psi_{t+s}\|$ by (2.8), this means that $\varphi_t - \psi_{t+s}$ is bounded. But $\psi$ is bounded, so $\varphi$ is also bounded.

Let $\rho_{0T}(f, g) := \sup_{0 \leq t \leq T} \|f - g\|$ $(f, g \in C_{0T})$ be the metric corresponding to uniform topology, and denote by $H_y^1$ the Cameron–Martin space of absolutely continuous functions starting at $y$ that possess square-integrable derivatives.

PROPOSITION 3.3. *The family $(Y^{s,y})$ satisfies a large deviations principle with good rate function*

$$(3.6) \quad I_{0T}^{s,y}(\varphi) = \begin{cases} \frac{1}{2}\int_0^T \|\dot{\varphi}_t - V(\varphi_t) + \Phi(\varphi_t - \psi_{t+s}(x_0))\|^2 \, dt, & \text{if } \varphi \in H_y^1, \\ \infty, & \text{otherwise }. \end{cases}$$

*More precisely, for any closed set $F \subset C_{0T}$ we have*

$$\limsup_{\varepsilon \to 0} \varepsilon \log \mathbb{P}(Y^{s,y} \in F) \leq -\inf_{\phi \in F} I_{0T}^{s,y}(\phi),$$

*and for any open set $G \subset C_{0T}$*

$$\liminf_{\varepsilon \to 0} \varepsilon \log \mathbb{P}(Y^{s,y} \in G) \geq -\inf_{\phi \in G} I_{0T}^{s,y}(\phi).$$

PROOF. Let $a(t, y) := V(y) - \Phi(y - \psi_t)$, and denote by $F$ the function that maps a path $g \in C_{0T}$ to the solution $f$ of the ODE

$$f_t = x_0 + \int_0^t a(s, f_s) \, ds + g_t, \qquad 0 \leq t \leq T.$$

Fix $g \in C_{0T}$, and let $R > 0$ such that the deterministic trajectory $\psi(x_0)$ as well as $f = F(g)$ stay in $B_R(0)$ up to time $T$. Note that nonexplosion of $f$ is guaranteed by dissipativity of $a$, much like in (3.5). Now observe that $a$ is locally Lipschitz with constant $2K_{2R}$ on $B_R(0)$, uniformly w.r.t. $t \in [0, T]$. Thus, we have for $\tilde{g} \in C_{0T}$, $\tilde{f} = F(\tilde{g})$ such that $\tilde{f}$ does not leave $B_R(0)$ up to time $T$:

$$\|f_t - \tilde{f}_t\| \leq 2K_{2R} \int_0^t \|f_s - \tilde{f}_s\| \, ds + \|g_t - \tilde{g}_t\|,$$

and Gronwall's lemma yields

$$\rho_{0T}(f, \tilde{f}) \leq \rho_{0T}(g, \tilde{g}) e^{2K_{2R}T},$$

that is, $F$ is continuous. Indeed, the last inequality shows that we do not have to presume that $\tilde{f}$ stays in $B_R(0)$, but that this is granted whenever $\rho_{0T}(g, \tilde{g})$ is sufficiently small.

Since $F$ is continuous and $F(\sqrt{\varepsilon}W) = Y$, we may invoke Schilder's theorem and the contraction principle, to deduce that $Y$ satisfies a large deviations principle with rate function

$$I_{0T}(\varphi) = \inf\left\{\tfrac{1}{2} \int_0^T \|\dot{g}_t\|^2 \, dt : g \in H_y^1, F(g) = \varphi\right\}.$$

This proves the LDP for $(Y^{s,y})$.   $\square$

Notice that the rate function of $Y$ measures distances from the deterministic solution $\psi$ just as in the classical case without interaction, but the distance of $\varphi$ from $\psi$ is weighted by the interaction between the two paths.

By means of the rate function, one can associate to $Y^{s,y}$ two functions that determine the cost, respectively energy, of moving between points in the geometric landscape induced by the vector field $V$. For $t \geq 0$, the *cost function*

$$C^s(y, z, t) = \inf_{f \in C_{0t} : f_t = z} I_{0t}^{s,y}(f), \qquad y, z \in \mathbb{R}^d,$$

determines the asymptotic cost for the diffusion $Y^{s,y}$ to move from $y$ to $z$ in time $t$, and the *quasi-potential*

$$Q^s(y, z) = \inf_{t > 0} C^s(y, z, t)$$

describes its cost of going from $y$ to $z$ eventually.

### 3.2. *Large deviations principle for the self-stabilizing diffusion.*

We are now in a position to prove large deviations principles for $\xi$ and $X$ by showing that $\xi$ and $Y$ are close in the sense of large deviations.

THEOREM 3.4. *For any $\varepsilon > 0$ let $x_0^\varepsilon, \xi_0^\varepsilon \in \mathbb{R}^d$ that converge to some $x_0 \in \mathbb{R}^d$, respectively $y \in \mathbb{R}^d$, as $\varepsilon \to 0$. Denote by $X^\varepsilon$ the solution of (3.1) starting at $x_0^\varepsilon$. Let $s \geq 0$, and denote by $\xi^\varepsilon$ the solution of (3.3) starting in $\xi_0^\varepsilon$ with time parameter $s$, that is,*

$$(3.7) \qquad \xi_t^\varepsilon = \xi_0^\varepsilon + \int_0^t V(\xi_u^\varepsilon) - b^{\varepsilon,x_0}(u+s,\xi_u^\varepsilon)\, du + \sqrt{\varepsilon} W_t, \qquad t \geq 0,$$

*where $b^{\varepsilon,x_0}(t,x) = \mathbb{E}[\Phi(x - X_t^\varepsilon)]$.*

*Then the diffusions $(\xi^\varepsilon)_{\varepsilon>0}$ satisfy on any time interval $[0,T]$ a large deviations principle with good rate function (3.6).*

PROOF. We shall show that $\xi := \xi^\varepsilon$ is exponentially equivalent to $Y := Y^{s,y}$ as defined by (3.4), which has the desired rate function; that is, we prove that for any $\delta > 0$ we have

$$(3.8) \qquad \limsup_{\varepsilon \to 0} \varepsilon \log \mathbb{P}(\rho_{0T}(\xi, Y) \geq \delta) = -\infty.$$

Without loss of generality, we may choose $R > 0$ such that $x_0^\varepsilon, y \in B_R(0)$ and that $\psi_t(x_0)$ does not leave $B_R(0)$ up to time $s+T$, and denote by $\sigma_R$ the first time at which $\xi$ or $Y$ exits from $B_R(0)$. Then for $t \leq \sigma_R$

$$(3.9) \qquad \begin{aligned} \|\xi_t - Y_t\| &\leq \|\xi_0 - y\| + \int_0^t \|V(\xi_u) - V(Y_u)\|\, du \\ &\quad + \int_0^t \|b^{\varepsilon,x_0^\varepsilon}(u+s,\xi_u) - \Phi(Y_u - \psi_{u+s}(x_0))\|\, du. \end{aligned}$$

The first integral satisfies

$$\int_0^t \|V(\xi_u) - V(Y_u)\|\, du \leq K_R \int_0^t \|\xi_u - Y_u\|\, du, t \leq \sigma_R,$$

due to the local Lipschitz assumption. Let us decompose the second integral. We have

$$\begin{aligned} \|b^{\varepsilon,x_0^\varepsilon}&(u+s,\xi_u) - \Phi(Y_u - \psi_{u+s}(x_0))\| \\ &\leq \|b^{\varepsilon,x_0^\varepsilon}(u+s,\xi_u) - \Phi(\xi_u - \psi_{u+s}(x_0^\varepsilon))\| \\ &\quad + \|\Phi(\xi_u - \psi_{u+s}(x_0^\varepsilon)) - \Phi(\xi_u - \psi_{u+s}(x_0))\| \\ &\quad + \|\Phi(\xi_u - \psi_{u+s}(x_0)) - \Phi(Y_u - \psi_{u+s}(x_0))\|. \end{aligned}$$

Bounds for the second and third term in this decomposition are easily derived. The last one is seen to be bounded by $K_{2R}\|\xi_u - Y_u\|$, since $\xi, Y$ as well as $\psi$ are in $B_R(0)$ before time $\sigma_R \wedge T$. For the second term we also use the Lipschitz condition to deduce that

$$\|\Phi(\xi_u - \psi_{u+s}(x_0^\varepsilon)) - \Phi(\xi_u - \psi_{u+s}(x_0))\| \leq K_{2R}\|\psi_{u+s}(x_0^\varepsilon) - \psi_{u+s}(x_0)\|.$$

As a consequence of the flow property for $\psi$ this bound approaches 0 as $\varepsilon \to 0$ uniformly w.r.t. $u \in [0, T]$.

By combining these bounds and applying Gronwall's lemma, we find that

$$
\begin{aligned}
& \|\xi_t - Y_t\| \\
(3.10) \quad & \leq \exp\{2K_{2R}t\} \bigg( \|\xi_0 - y\| + K_{2R} \int_0^t \|\psi_{u+s}(x_0^\varepsilon) - \psi_{u+s}(x_0)\| \, du \\
& \qquad\qquad + \int_0^t \|b^{\varepsilon, x_0^\varepsilon}(u+s, \xi_u) - \Phi(\xi_u - \psi_{u+s}(x_0^\varepsilon))\| \, du \bigg)
\end{aligned}
$$

for $t \leq \sigma_R$. Since $\xi$ is bounded before $\sigma_R$ the r.h.s. of this inequality tends to zero by Corollary 3.2.

The exponential equivalence follows from the LDP for $Y$ as follows. Fix $\delta > 0$, and choose $\varepsilon_0 > 0$ such that the r.h.s. of (3.10) is smaller than $\delta$ for $\varepsilon \leq \varepsilon_0$. Then $\|\xi_t - Y_t\| > \delta$ implies that at least one of $\xi_t$ or $Y_t$ is not in $B_R(0)$, and if $\xi_t \notin B_R(0)$, then $Y_t \notin B_{R/2}(0)$ if $\delta$ is small enough. Thus we can bound the distance of $\xi$ and $Y$ by an exit probability of $Y$. For $l > 0$ let $\tau_l$ denote the diffusion $Y$'s time of first exit from $B_l(0)$. Then, by Proposition 3.3,

$$
\begin{aligned}
& \limsup_{\varepsilon \to 0} \varepsilon \log \mathbb{P}(\rho_{0T}(\xi, Y) > \delta) \\
(3.11) \quad & \leq \limsup_{\varepsilon \to 0} \varepsilon \log \mathbb{P}(\tau_{R/2} \leq T) \\
& \leq -\inf\left\{ C^s(y, z, t) : |z| \geq \frac{R}{2}, 0 \leq t \leq T \right\}.
\end{aligned}
$$

The latter expression approaches $-\infty$ as $R \to \infty$.  $\square$

Theorem 3.4 allows us to deduce two important corollaries. A particular choice of parameters yields an LDP for $X$, and the $\varepsilon$-dependence of the initial conditions permits us to conclude that the LDP holds uniformly on compact subsets, a fact that is crucial for the proof of an exit law in the following section. The arguments can be found in [7].

Let $\mathbb{P}_{x_0}(X \in \cdot)$ denote the law of the diffusion $X$ starting at $x_0 \in \mathbb{R}^d$.

COROLLARY 3.5.   *Let $L \subset \mathbb{R}^d$ be a compact set. For any closed set $F \subset C_{0T}$ we have*

$$
\limsup_{\varepsilon \to 0} \varepsilon \log \sup_{x_0 \in L} \mathbb{P}_{x_0}(X \in F) \leq -\inf_{x_0 \in L} \inf_{\phi \in F} I_{0T}^{0, x_0}(\phi),
$$

*and for any open set $G \subset C_{0T}$*

$$
\liminf_{\varepsilon \to 0} \varepsilon \log \inf_{x_0 \in L} \mathbb{P}_{x_0}(X \in G) \geq -\sup_{x_0 \in L} \inf_{\phi \in G} I_{0T}^{0, x_0}(\phi).
$$

PROOF. Choosing $x_0^\varepsilon = \xi_0^\varepsilon$ and $s = 0$ implies $\xi^\varepsilon = X^\varepsilon$ in Theorem 3.4, which shows that $X$ satisfies an LDP with rate function $I_{0T}^{0,x_0}$. Furthermore, this LDP allows for $\varepsilon$-dependent initial conditions. This implies the uniformity of the LDP, as pointed out in the proofs of Theorem 5.6.12 and Corollary 5.6.15 in [7]. Indeed, the $\varepsilon$-dependence yields for all $x_0 \in \mathbb{R}^d$

$$\limsup_{\varepsilon \to 0, y \to x_0} \varepsilon \log \mathbb{P}_y(X \in F) \leq -\inf_{\phi \in F} I_{0T}^{0,x_0}(\phi),$$

for otherwise one could find sequences $\varepsilon_n > 0$ and $y_n \in \mathbb{R}^d$ such that $\varepsilon_n \to 0$, $y_n \to x_0$ and

$$\limsup_{n \to \infty} \varepsilon_n \log \mathbb{P}_{y_n}(X \in F) > -\inf_{\phi \in F} I_{0T}^{0,x_0}(\phi).$$

But this contradicts the LDP.

Now the uniformity of the upper large deviations bound follows exactly as demonstrated in the proof of Corollary 5.6.15 in [7]. The lower bound is treated similarly. □

The next corollary is just a consequence of the $\varepsilon$-dependent initial conditions in the LDP for $\xi$.

COROLLARY 3.6. *Let* $L \subset \mathbb{R}^d$ *be a compact set. For any closed set* $F \subset C_{0T}$ *we have*

$$\limsup_{\varepsilon \to 0} \varepsilon \log \sup_{x_0 \in L} \mathbb{P}(\xi^{s,x_0} \in F) \leq -\inf_{x_0 \in L} \inf_{\phi \in F} I_{0T}^{s,x_0}(\phi),$$

*and for any open set* $G \subset C_{0T}$

$$\liminf_{\varepsilon \to 0} \varepsilon \log \inf_{x_0 \in L} \mathbb{P}(\xi^{s,x_0} \in G) \geq -\sup_{x_0 \in L} \inf_{\phi \in G} I_{0T}^{s,x_0}(\phi).$$

3.3. *Exponential approximations under stability assumptions.* The aim of this subsection is to exploit the fact that the inhomogeneity of the diffusion $Y^{s,y}$ is weak in the sense that its drift depends on time only through $\psi_{t+s}(x_0)$. If the dynamical system $\dot{\psi} = V(\psi)$ admits an asymptotically stable fixed point $x_{\text{stable}}$ that attracts $x_0$, then the drift of $Y^{s,y}$ becomes almost autonomous for large times, which in turn may be used to estimate large deviations probabilities for $\xi^{s,y}$. We make the following assumption. It will also be in force in Section 4, where it will keep us from formulating results on exits from domains with boundaries containing critical points of $DV$, in particular saddle points in the potential case.

ASSUMPTION 3.7.

(i) Stability: there exists a stable equilibrium point $x_{\text{stable}} \in \mathbb{R}^d$ of the dynamical system

$$\dot{\psi} = V(\psi).$$

(ii) Convexity: the geometry induced by the vector field $V$ is convex, that is, the condition (2.6) for $V$ holds globally:

$$(3.12) \qquad \langle h, DV(x)h \rangle \leq -K_V$$

for $h \in \mathbb{R}^d$ s.t. $\|h\| = 1$ and $x \in \mathbb{R}^d$.

Under this assumption it is natural to consider the limiting time-homogeneous diffusion $Y^{\infty,y}$ defined by

$$(3.13) \quad dY_t^\infty = V(Y_t^\infty)\,dt - \Phi(Y_t^\infty - x_{\text{stable}})\,dt + \sqrt{\varepsilon}\,dW_t, \qquad Y_0^\infty = y.$$

LEMMA 3.8.    *Let $L \subset \mathbb{R}^d$ be compact, and assume that $x_{\text{stable}}$ attracts all $y \in L$, that is,*

$$\lim_{t \to \infty} \psi_t(y) = x_{\text{stable}} \qquad \forall y \in L.$$

*Then $Y^{\infty,y}$ is an exponentially good approximation of $Y^{s,y}$, that is, for any $\delta > 0$ we have*

$$\lim_{r \to \infty} \limsup_{\varepsilon \to 0} \varepsilon \log \sup_{y \in L, s \geq r} \mathbb{P}(\rho_{0T}(Y^{s,y}, Y^{\infty,y}) \geq \delta) = -\infty.$$

PROOF.    We have

$$\|Y_t^{s,y} - Y_t^{\infty,y}\| \leq \int_0^t \|V(Y_u^{s,y}) - V(Y_u^{\infty,y})\|\,du$$

$$+ \int_0^t \|\Phi(Y_u^{s,y} - \psi_{s+u}(y)) - \Phi(Y_u^{\infty,y} - x_{\text{stable}})\|\,du.$$

Let $\sigma_R^{s,y}$ be the first time at which $Y^{s,y}$ or $Y^{\infty,y}$ exits from $B_R(0)$. For $t \leq \sigma_R^{s,y}$, we may use the Lipschitz property of $\Phi$ and $V$, to find a constant $c_R > 0$ s.t.

$$\|Y_t^{s,y} - Y_t^\infty\| \leq c_R \int_0^t \|Y_u^{s,y} - Y_u^\infty\|\,du + c_R T \rho_{0T}(\psi_{s+\cdot}(y), x_{\text{stable}}).$$

By assumption the second term converges to 0 as $s \to \infty$, uniformly with respect to $y \in L$ since the flow is continuous with respect to the initial data. Hence, by Gronwall's lemma there exists some $r = r(R, \delta) > 0$ such that for $s \geq r$

$$\sup_{y \in L} \sup_{0 \leq t \leq \sigma_R^{s,y}} \|Y_t^{s,y} - Y_t^\infty\| < \delta/2.$$

We deduce that

$$\mathbb{P}(\rho_{0T}(Y^{s,y}, Y^{\infty}) \geq \delta/2) \leq \mathbb{P}(\tau_{R/2}^{y} \leq T) \qquad \forall s \geq r, y \in L,$$

where for $l > 0$ $\tau_l^y$ denotes the first exit time of $Y^{\infty,y}$ from $B_l(0)$. Sending $r, R \to \infty$ and appealing to the uniform LDP for $Y^{\infty,y}$ finishes the proof, much as the proof of Theorem 3.4. $\quad\square$

This exponential closeness of $Y^{\infty,y}$ and $Y^{s,y}$ carries over to $\xi^{s,y}$ under the aforementioned stability and convexity assumption, which enables us to sharpen the exponential equivalence proved in Theorem 3.4. In order to establish this improvement, we need a preparatory lemma that strengthens Corollary 3.2 to uniform convergence over the whole time axis. This uniformity is of crucial importance for the proof of an exit law in the next section and depends substantially on the strong convexity assumption (3.12).

LEMMA 3.9.  *We have*

$$\lim_{\varepsilon \to 0} b^{\varepsilon, x_0}(t, x) = \Phi(x - \psi_t(x_0)),$$

*uniformly w.r.t. $t \geq 0$ and w.r.t. $x$ and $x_0$ on compact subsets of $\mathbb{R}^d$.*

PROOF.  Let $f(t) := \mathbb{E}(\|Z_t\|^2)$, where $Z_t = X_t - \psi_t(x_0)$. In the proof of Lemma 3.1 we have seen that

$$f'(t) \leq 2\mathbb{E}[\langle Z_t, V(Z_t + \psi_t(x_0)) - V(\psi_t(x_0))\rangle] + \varepsilon d$$

$$\leq -2K_V \mathbb{E}(\|Z_t\|^2) + \varepsilon d = -2K_V f(t) + \varepsilon d.$$

This means that $\{t \geq 0 : f'(t) < 0\} \supset \{t \geq 0 : f(t) > \frac{\varepsilon d}{2K_V}\}$. Recalling that $f(0) = 0$, this allows us to conclude that $f$ is bounded by $\frac{\varepsilon d}{2K_V}$. Now an appeal to the proof of Corollary 3.2 finishes the argument. $\quad\square$

PROPOSITION 3.10.  *Let $L \subset \mathbb{R}^d$ be compact, and assume that $x_{\text{stable}}$ attracts all $y \in L$. Then $Y^{\infty,y}$ is an exponentially good approximation of $\xi^{s,y}$, that is, for any $\delta > 0$ we have*

$$\lim_{r \to \infty} \limsup_{\varepsilon \to 0} \varepsilon \log \sup_{y \in L, s \geq r} \mathbb{P}(\rho_{0T}(\xi^{s,y}, Y^{\infty,y}) \geq \delta) = -\infty.$$

PROOF.  Recall the proof of Theorem 3.4. For $y \in L$ and $s \geq 0$ we have

$$\|\xi_t^{s,y} - Y_t^{s,y}\|$$

$$(3.14) \qquad \leq \exp\{2K_{2R}t\} \int_0^t \|b^{\varepsilon,x_0}(u + s, \xi_u) - \Phi(\xi_u - \psi_{u+s}(x_0))\| \, du$$

for $t \le \sigma_R^{y,s}$, which denotes the first time that $\xi_t^{s,y}$ or $Y^{s,y}$ exits from $B_R(0)$. By Lemma 3.9, the integrand on the r.h.s. converges to zero as $\varepsilon \to 0$, uniformly w.r.t. $s \ge 0$. Therefore, if we fix $\delta > 0$, we may choose $R = R(\delta)$ sufficiently large and $\varepsilon_0 > 0$ such that for $\varepsilon \le \varepsilon_0$, and all $s \ge 0$

$$
\mathbb{P}(\rho_{0T}(\xi^{s,y}, Y^{s,y}) > \delta) \le \mathbb{P}(\tau_{R/2}^{s,y} < T)
$$
$$
\le \mathbb{P}(\tau_{R/4}^{\infty,y} < T) + \mathbb{P}(\rho_{0T}(Y^{\infty,y}, Y^{s,y}) > R/4),
$$

where for $l > 0$, $0 \le s \le \infty$, $\tau_l^{s,y}$ denotes the first exit time of the diffusion $Y^{s,y}$ from the ball $B_l(0)$. By the uniform LDP for $Y^{\infty,y}$ and Lemma 3.8 the assertion follows. $\square$

## 4. The exit problem.

As a consequence of the large deviations principle, the trajectories of the self-stabilizing diffusion are attracted to the deterministic dynamical system $\dot{\psi} = V(\psi)$ as noise tends to 0. The probabilities of deviating from $\psi$ are exponentially small in $\varepsilon$, and the diffusion will certainly exit from a domain within a certain time interval if the deterministic path $\psi$ exits. The problem of diffusion exit involves an analysis for the rare event that the diffusion leaves the domain although the deterministic path stays inside, that is, it is concerned with an exit which is triggered by noise only. Clearly, the time of such an exit should increase as the noise intensity tends to zero. In this section we shall derive the precise large deviations asymptotics of such exit times, that is, we shall give an analogue of the well-known Kramers–Eyring law for time-homogeneous diffusions.

Let us briefly recall this law, a detailed presentation of which may be found in Section 5.7 of [7]. For further classical results about the exit problem we refer to [5, 6, 8] and [17].

A Brownian particle of intensity $\varepsilon$ that wanders in a geometric landscape given by a potential $U$ is mathematically described by the classical time-homogeneous SDE

$$
dZ_t^\varepsilon = -\nabla U(Z_t^\varepsilon)\, dt + \sqrt{\varepsilon}\, dW_t, \qquad Z_0^\varepsilon = x_0 \in \mathbb{R}^d.
$$

If $x^*$ is a stable fixed point of the system $\dot{x} = -\nabla U(x)$ that attracts the initial condition $x_0$ and $\tau^\varepsilon$ denotes the exit time from the domain of attraction of $x^*$, then the asymptotics of $\tau^\varepsilon$ is described by the following two relations:

$$
(4.1) \qquad \lim_{\varepsilon \to 0} \varepsilon \log \mathbb{E}(\tau^\varepsilon) = \bar{U},
$$

$$
(4.2) \qquad \lim_{\varepsilon \to 0} \mathbb{P}(e^{(\bar{U}-\delta)/\varepsilon} < \tau^\varepsilon < e^{(\bar{U}+\delta)/\varepsilon}) = 1 \qquad \forall \delta > 0.
$$

Here $\bar{U}$ denotes the energy required to exit from the domain of attraction of $x^*$. This law may roughly be paraphrased by saying that $\tau^\varepsilon$ behaves like $\exp \frac{\bar{U}}{\varepsilon}$ as $\varepsilon \to 0$.

Let us now return to the self-stabilizing diffusion $X^\varepsilon$, defined by (3.1). Intuitively, exit times should increase compared to the classical case due to self-stabilization and the inertia it entails. We shall show that this is indeed the case, and prove a synonym of (4.1) and (4.2) for the self-stabilizing diffusion. Our approach follows the presentation in [7].

Let $D$ be an open bounded domain in $\mathbb{R}^d$ in which $X^\varepsilon$ starts, that is, $x_0 \in D$, and denote by

$$\tau_D^\varepsilon = \inf\{t > 0 \colon X_t^\varepsilon \in \partial D\}$$

the first exit time from $D$. We make the following stability assumptions about $D$.

ASSUMPTION 4.1.

(i) The unique equilibrium point in $D$ of the dynamical system

$$(4.3) \qquad \dot{\psi}_t = V(\psi_t)$$

is stable and given by $x_{\text{stable}} \in D$. As before, $\psi_t(x_0)$ denotes the solution starting at $x_0$. We assume that $\lim_{t \to \infty} \psi_t(x_0) = x_{\text{stable}}$.

(ii) The solutions of

$$(4.4) \qquad \dot{\phi}_t = V(\phi_t) - \Phi(\phi_t - x_{\text{stable}})$$

satisfy

$$\phi_0 \in D \quad \Longrightarrow \quad \phi_t \in D \qquad \forall t > 0,$$

$$\phi_0 \in \overline{D} \quad \Longrightarrow \quad \lim_{t \to \infty} \phi_t = x_{\text{stable}}.$$

The description of the exponential rate for the exit time of Itô diffusions with homogeneous coefficients was first proved by Freidlin and Wentzell via an exploitation of the strong Markov property. The self-stabilizing diffusion $X^\varepsilon$ is also Markovian, but it is inhomogeneous, which makes a direct application of the Markov property difficult. However, the inhomogeneity is weak under the stability Assumption 4.1. It implies that the law of $X_t^\varepsilon$ converges as time tends to infinity, and large deviations probabilities for $X^\varepsilon$ may be approximated by those of $Y^\infty$ in the sense of Proposition 3.10. Since $Y^\infty$ is defined in terms of an autonomous SDE, its exit behavior is accessible through classical results. The rate function that describes the LDP for $Y^\infty$ is given by

$$(4.5) \quad I_{0T}^{\infty,y}(\varphi) = \begin{cases} \dfrac{1}{2} \displaystyle\int_0^T \|\dot{\varphi}_t - V(\varphi_t) + \Phi(\varphi_t - x_{\text{stable}})\|^2 \, dt, & \text{if } \varphi \in H_y^1, \\ \infty, & \text{otherwise.} \end{cases}$$

The corresponding cost function and quasi-potential are defined in an obvious way and denoted by $C^\infty$ and $Q^\infty$, respectively. The minimal energy required to connect the stable equilibrium point $x_{\text{stable}}$ to the boundary of the domain is assumed to be finite, that is,

$$\overline{Q}_\infty := \inf_{z \in \partial D} Q^\infty(x_{\text{stable}}, z) < \infty.$$

The following two theorems state our main result about the exponential rate of the exit time and the exit location.

THEOREM 4.2.    *For all $x_0 \in D$ and all $\eta > 0$, we have*

$$(4.6) \qquad \limsup_{\varepsilon \to 0} \varepsilon \log\{1 - \mathbb{P}_{x_0}(e^{(\overline{Q}_\infty - \eta)/\varepsilon} < \tau_D^\varepsilon < e^{(\overline{Q}_\infty + \eta)/\varepsilon})\} \le -\eta/2$$

*and*

$$(4.7) \qquad \lim_{\varepsilon \to 0} \varepsilon \log \mathbb{E}_{x_0}(\tau_D^\varepsilon) = \overline{Q}_\infty.$$

THEOREM 4.3.    *If $N \subset \partial D$ is a closed set satisfying*

$$\inf_{z \in N} Q^\infty(x_{\text{stable}}, z) > \overline{Q}_\infty,$$

*then it does not see the exit point: for any $x_0 \in D$*

$$\lim_{\varepsilon \to 0} \mathbb{P}_{x_0}(X_{\tau_D}^\varepsilon \in N) = 0.$$

The rest of this section is devoted to the proof of these two theorems. In the subsequent section, these results are illustrated by examples which show that the attraction part of the drift term in a diffusion may completely change the behavior of the paths, that is, the self-stabilizing diffusion stays in the domain for a longer time than the classical one, and it typically exits at a different place.

4.1. *Enlargement of the domain.*  The self-stabilizing diffusion lives in the open, bounded domain $D$ which is assumed to fulfill the previously stated stability conditions. In order to derive upper and lower bounds of exit probabilities, we need to construct an enlargement of $D$ that still enjoys the stability properties of Assumption 4.1(ii). This is possible because the family of solutions to the dynamical system (4.4) defines a continuous flow.

For $\delta > 0$ we denote by $D^\delta := \{y \in \mathbb{R}^d : \text{dist}(y, D) < \delta\}$ the open $\delta$-neighborhood of $D$. The flow $\phi$ is continuous, hence uniformly continuous on $\overline{D}$ due to boundedness of $D$, and since the vector field is locally Lipschitz. Hence, if $\delta$ is small enough, the trajectories $\phi_t(y)$ converge to $x_{\text{stable}}$ for $y \in D^\delta$, that is, for each neighborhood $\mathcal{V} \subset D$ of $x_{\text{stable}}$ there exists some $T > 0$ such that for $y \in D^\delta$ we have $\phi_t(y) \in \mathcal{V}$ for all $t \ge T$. Moreover, the joint continuity of the flow implies that, if we fix $c > 0$, we may choose $\delta = \delta(c) > 0$ such that

$$\sup\{\text{dist}(\phi_t(y), D) : t \in [0, T], y \in D^\delta\} < c.$$

Let

$$\mathcal{O}^\delta = \left\{ y \in \mathbb{R}^d : \sup_{t \in [0,T]} \mathrm{dist}(\phi_t(y), D) < c, \phi_T(y) \in \mathcal{V} \right\}.$$

Then $\mathcal{O}^\delta$ is a bounded open set which contains $D^\delta$ and satisfies Assumption 4.1(ii). Indeed, if $\delta$ is small enough, the boundary of $\mathcal{O}^\delta$ is not a characteristic boundary, and $\bigcap_{\delta>0} \mathcal{O}^\delta = D$.

4.2. *Proof of the upper bound for the exit time.* For the proof of the two main results, we successively proceed in several steps and establish a series of preparatory estimates that shall be combined afterward. In this subsection, we concentrate on the upper bound for the exit time from $D$, and establish inequalities for the probability of exceeding this bound and for the mean exit time.

In the sequel, we denote by $\mathbb{P}_{s,y}$ the law of the diffusion $\xi^{s,y}$, defined by (3.3). Recall that by the results of the previous section, $\xi^{s,y}$ satisfies a large deviations principle with rate function $I^{s,y}$. The following continuity property of the associated cost function is the analogue of Lemma 5.7.8 in [7] for this inhomogeneous diffusion. The proof is omitted.

LEMMA 4.4. *For any $\delta > 0$ and $s \in [0, \infty)$, there exists $\varrho > 0$ such that*

$$(4.8) \qquad \sup_{x,y \in B_\varrho(x_{\text{stable}})} \inf_{t \in [0,1]} C^s(x, y, t) < \delta$$

*and*

$$(4.9) \qquad \sup_{(x,y) \in \Gamma} \inf_{t \in [0,1]} C^s(x, y, t) < \delta,$$

*where $\Gamma = \{(x, y) : \inf_{z \in \partial D}(\|y - z\| + \|x - z\|) \leq \varrho\}$.*

Let us now present two preliminary lemmas on exit times of $\xi^{s,y}$. In slight abuse of notation, we denote exit times of $\xi^{s,y}$ also by $\tau_D^\varepsilon$, which could formally be justified by assuming to look solely at the coordinate process on path space and switching between measures instead of processes. On the other hand, this notation is convenient when having in mind that $\xi^{s,y}$ describes the law of $X^\varepsilon$ restarted at time $s$, and that $X^\varepsilon$ may be recovered from $\xi^{s,y}$ for certain parameters.

LEMMA 4.5. *For any $\eta > 0$ and $\varrho > 0$ small enough, there exist $T_0 > 0$, $s_0 > 0$ and $\varepsilon_0 > 0$ such that*

$$\inf_{y \in B_\varrho(x_{\text{stable}})} \mathbb{P}_{s,y}(\tau_D^\varepsilon \leq T_0) \geq e^{-(\overline{Q}_\infty + \eta)/\varepsilon} \qquad \text{for all } \varepsilon \leq \varepsilon_0 \text{ and } s \geq s_0.$$

PROOF. Let $\varrho$ be given according to Lemma 4.4. The corresponding result for the time-homogeneous diffusion $Y^{\infty,y}$ is well known (see [7], Lemma 5.7.18), and will be carried over to $\xi^{s,y}$ using the exponential approximation of Proposition 3.10. Let $\mathbb{P}_{\infty,y}$ denote the law of $Y^{\infty,y}$. The drift of $Y^{\infty,y}$ is locally Lipschitz by the assumptions on $V$ and $\Phi$, and we may assume w.l.o.g. that it is even globally Lipschitz. Otherwise we change the drift outside a large domain containing $D$.

If $\delta > 0$ is small enough such that the enlarged domain $\mathcal{O}^\delta$ satisfies Assumption 4.1(ii), Lemma 5.7.18 in [7] implies the existence of $\varepsilon_1$ and $T_0$ such that

$$(4.10) \qquad \inf_{y \in B_\varrho(x_{\text{stable}})} \mathbb{P}_{\infty,y}(\tau_{\mathcal{O}^\delta}^\varepsilon \leq T_0) \geq e^{-(\overline{Q}_\infty^\delta + \eta/3)/\varepsilon} \qquad \text{for all } \varepsilon \leq \varepsilon_1.$$

Here $\overline{Q}_\infty^\delta$ denotes the minimal energy

$$\overline{Q}_\infty^\delta = \inf_{z \in \partial \mathcal{O}^\delta} Q^\infty(x_{\text{stable}}, z).$$

The continuity of the cost function carries over to the quasi-potential, that is, there exists some $\delta_0 > 0$ such that $|\overline{Q}_\infty^\delta - \overline{Q}_\infty| \leq \eta/3$ for $\delta \leq \delta_0$.

Now let us link the exit probabilities of $Y^{\infty,y}$ and $\xi^{s,y}$. We have for $s \geq 0$

$$\begin{aligned}
& \mathbb{P}_{s,y}(\tau_D^\varepsilon \leq T_0) \\
(4.11) \quad & \geq \mathbb{P}(\{\xi^{s,y} \text{ exits from } D \text{ before } T_0\} \cap \{\rho_{0,T_0}(\xi^{s,y}, Y^{\infty,y}) \leq \delta\}) \\
& \geq \mathbb{P}_{\infty,y}(\tau_{D^\delta}^\varepsilon \leq T_0) - \mathbb{P}(\rho_{0,T_0}(\xi^{s,y}, Y^\infty) \geq \delta).
\end{aligned}$$

Moreover, by the exponential approximation we may find $\varepsilon_2 > 0$ and $s_0 > 0$ such that

$$\sup_{y \in B_\varrho(x_{\text{stable}})} \mathbb{P}(\rho_{0,T_0}(\xi^{s,y}, Y^\infty) \geq \delta) \leq e^{-(\overline{Q}_\infty^\delta + \eta/2)/\varepsilon} \qquad \forall s \geq s_0, \varepsilon \leq \varepsilon_2.$$

Since $D^\delta \subset \mathcal{O}^\delta$, we deduce that for $\varepsilon \leq \varepsilon_0 = \varepsilon_1 \wedge \varepsilon_2$ and $s \geq s_0$

$$\inf_{y \in B_\varrho(x_{\text{stable}})} \mathbb{P}_{s,y}(\tau_D^\varepsilon \leq T_0) \geq e^{-(\overline{Q}_\infty^\delta + \eta/3)/\varepsilon} - e^{-(\overline{Q}_\infty^\delta + \eta/2)/\varepsilon} \geq e^{-(\overline{Q}_\infty^\delta + \eta)/\varepsilon}. \quad \square$$

By similar arguments, we prove the exponential smallness of the probability of too long exit times. Let $\Sigma_\varrho = \inf\{t \geq 0 : \xi_t^{s,y} \in B_\varrho(x_{\text{stable}}) \cup \partial D\}$, where $\varrho$ is small enough such that $B_\varrho(x_{\text{stable}})$ is contained in the domain $D$.

LEMMA 4.6. *For any $\varrho > 0$ sufficiently small and for any $K > 0$ there exist $\varepsilon_0 > 0$, $T_1 > 0$ and $r > 0$ such that*

$$\sup_{y \in D, s \geq r} \mathbb{P}_{s,y}(\Sigma_\varrho > t) \leq e^{-K/\varepsilon} \qquad \forall t \geq T_1.$$

PROOF. As before, we use the fact that a similar result is already known for $Y^{\infty,y}$. For $\delta > 0$ small enough, let

$$\Sigma_\varrho^\delta = \inf\{t \geq 0 : Y_t^\infty \in B_{\varrho-\delta}(x_{\text{stable}}) \cup \partial\mathcal{O}^\delta\}.$$

By Lemma 5.7.19 in [7], there exist $T_1 > 0$ and $\varepsilon_1 > 0$ such that

$$\sup_{y \in D} \mathbb{P}_{\infty,y}(\Sigma_\varrho^\delta > t) \leq e^{-K/\varepsilon} \qquad \forall t \geq T_1; \varepsilon \leq \varepsilon_1.$$

Now the assertion follows from

$$\sup_{y \in D} \mathbb{P}_{s,y}(\Sigma_\varrho > T_1) \leq \sup_{y \in D} \mathbb{P}_{\infty,y}(\Sigma_\varrho^\delta > T_1) + \sup_{y \in D} \mathbb{P}(\rho_{0,T_1}(\xi^{s,y}, Y^{\infty,y}) > \delta),$$

since the last term is exponentially negligible by Proposition 3.10. $\square$

The previous two lemmas contain the essential large deviations bounds required for the proof of the following upper bound for the exit time of $X^\varepsilon$.

PROPOSITION 4.7. *For all $x_0 \in D$ and $\eta > 0$ we have*

$$\limsup_{\varepsilon \to 0} \varepsilon \log \mathbb{P}_{x_0}(\tau_D^\varepsilon \geq e^{(\overline{Q}_\infty + \eta)/\varepsilon}) \leq -\eta/2 \tag{4.12}$$

*and*

$$\limsup_{\varepsilon \to 0} \varepsilon \log \mathbb{E}_{x_0}[\tau_D^\varepsilon] \leq \overline{Q}_\infty. \tag{4.13}$$

PROOF. The proof consists of a careful modification of the arguments used in Theorem 5.7.11 in [7]. By Lemmas 4.5 and 4.6, there exist $\widetilde{T} = T_0 + T_1 > 0$, $\varepsilon_0 > 0$ and $r_0 > 0$ such that for $T \geq \widetilde{T}$, $\varepsilon \leq \varepsilon_0$ and $r \geq r_0$ we have

$$
\begin{aligned}
q_T^r := \inf_{y \in D} \mathbb{P}_{r,y}(\tau_D^\varepsilon \leq T) & \\
\geq \inf_{y \in D} \mathbb{P}_{r,y}&(\Sigma_\varrho \leq T_1) \inf_{y \in B_\varrho(x_{\text{stable}}), s \geq r} \mathbb{P}_{s,y}(\tau_D^\varepsilon \leq T_0) \\
\geq \exp&\left\{-\frac{\overline{Q}_\infty + \eta/2}{\varepsilon}\right\} =: q_T^\infty.
\end{aligned}
\tag{4.14}
$$

Moreover, by the Markov property of $\xi^{s,y}$, we see that for $k \in \mathbb{N}$

$$
\begin{aligned}
\mathbb{P}_{x_0}(\tau_D^\varepsilon > 2(k+1)T) &= [1 - \mathbb{P}_{x_0}(\tau_D^\varepsilon \leq 2(k+1)T | \tau_D^\varepsilon > 2kT)] \\
&\quad \times \mathbb{P}_{x_0}(\tau_D^\varepsilon > 2kT) \\
&\leq \left[1 - \inf_{y \in D} \mathbb{P}_{2kT,y}(\tau_D^\varepsilon \leq 2T)\right] \mathbb{P}_{x_0}(\tau_D^\varepsilon > 2kT) \\
&\leq (1 - q_{2T}^{2kT}) \mathbb{P}_{x_0}(\tau_D^\varepsilon > 2kT),
\end{aligned}
$$

which by induction yields

$$(4.15) \qquad \mathbb{P}_{x_0}(\tau_D^\varepsilon > 2kT) \leq \prod_{i=0}^{k-1}(1 - q_{2T}^{2iT}).$$

Let us estimate each term of the product separately. We have

$$1 - q_{2T}^{2iT} = \sup_{y \in D} \mathbb{P}_{2iT,y}(\tau_D^\varepsilon > 2T)$$

$$\leq \sup_{y \in D} \mathbb{P}_{2iT,y}(\tau_D^\varepsilon > T) \sup_{y \in D} \mathbb{P}_{(2i+1)T,y}(\tau_D^\varepsilon > T)$$

$$\leq \sup_{y \in D} \mathbb{P}_{(2i+1)T,y}(\tau_D^\varepsilon > T).$$

By choosing $T$ large enough, we may replace the product in (4.15) by a power. Indeed, for $T > \max(\widetilde{T}, r_0)$ we have $(2i+1)T \geq r_0$ for all $i \in \mathbb{N}$, which by (4.14) results in the uniform upper bound

$$1 - q_{2T}^{2iT} \leq 1 - q_T^{(2i+1)T} \leq 1 - q_T^\infty.$$

By plugging this into (4.15), we obtain a "geometric" upper bound for the expected exit time, namely

$$\mathbb{E}_{x_0}[\tau_D^\varepsilon] \leq 2T \left[ 1 + \sum_{k=1}^\infty \sup_{y \in D} \mathbb{P}_{x_0}(\tau_D^\varepsilon \geq 2kT) \right]$$

$$\leq 2T \left[ 1 + \sum_{k=1}^\infty \prod_{i=0}^{k-1}(1 - q_{2T}^{2iT}) \right]$$

$$\leq 2T \left[ 1 + \sum_{k=1}^\infty (1 - q_T^\infty)^k \right] = \frac{2T}{q_T^\infty}.$$

This proves the claimed asymptotics of the expected exit time. Furthermore, an application of Chebyshev's inequality shows that

$$\mathbb{P}_{x_0}(\tau_D^\varepsilon \geq e^{(\overline{Q}_\infty + \eta)/\varepsilon}) \leq \frac{\mathbb{E}_{x_0}[\tau_D^\varepsilon]}{e^{(\overline{Q}_\infty + \eta)/\varepsilon}} \leq 2T \frac{e^{-(\overline{Q}_\infty + \eta)/\varepsilon}}{q_T^\infty} = 2Te^{-\eta/2\varepsilon},$$

which is the asserted upper bound of the exit probability. □

4.3. *Proof of the lower bound for the exit time.* In order to establish the lower bound of the exit time, we prove a preliminary lemma which estimates the probability to exit from the domain $D \setminus B_\varrho(x_{\text{stable}})$ at the boundary of $D$. This probability is seen to be exponentially small since the diffusion is attracted to the stable equilibrium point. Let us denote by $S_\varrho$ the boundary of $B_\varrho(x_{\text{stable}})$, and recall the definition of the stopping time $\Sigma_\varrho$.

Lemma 4.8. *For any closed set $N \subset \partial D$ and $\eta > 0$, there exist $\varepsilon_0 > 0$, $\varrho_0 > 0$ and $r_0 > 0$ such that*

$$\varepsilon \log \sup_{y \in S_{2\varrho}, s \geq r} \mathbb{P}(\xi_{\Sigma_\varrho^\infty}^{s,y} \in N) \leq -\inf_{z \in N} Q^\infty(x_{\text{stable}}, z) + \eta$$

*for all $\varepsilon \leq \varepsilon_0$, $r \geq r_0$ and $\varrho \leq \varrho_0$.*

Proof. For $\delta > 0$ we define a subset $\mathcal{S}^\delta$ of $D^\delta$ by setting

$$\mathcal{S}^\delta := D^\delta \setminus \{y \in \mathbb{R}^d : \text{dist}(y, N) < \delta\}.$$

Furthermore, let

$$\mathcal{N}^\delta := \partial \mathcal{S}^\delta \cap \{y \in \mathbb{R}^d : \text{dist}(y, N) \leq \delta\}.$$

$\mathcal{S}^\delta$ contains the stable equilibrium point $x_{\text{stable}}$, and as such it is unique in $\mathcal{S}^\delta$ if $\delta$ is small enough.

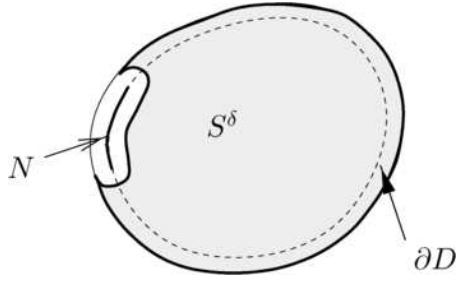

The proofs of Lemma 5.7.19 and Lemma 5.7.23 in [7] can be adapted to the domain $\mathcal{S}^\delta$, since an exit of the limiting diffusion $Y^\infty$ from the domain $\mathcal{O}^\delta$ defined in Section 4.1 always requires an exit from $\mathcal{S}^\delta$. Hence, there exist $\varepsilon_1 > 0$ and $\varrho_1 > 0$ such that

$$\varepsilon \log \sup_{y \in S_{2\varrho}} \mathbb{P}_{\infty, y}(Y_{\Sigma_\varrho^\delta}^\infty \in \mathcal{N}^\delta) \leq -\inf_{z \in \mathcal{N}^\delta} Q^\infty(x_{\text{stable}}, z) + \frac{\eta}{2}$$

for $\varepsilon \leq \varepsilon_1$ and $\varrho \leq \varrho_1$, where $\Sigma_\varrho^\delta$ denotes the first exit time from the domain $\mathcal{S}^\delta \setminus B_\varrho(x_{\text{stable}})$. By the continuity of the quasi-potential, we have

$$-\inf_{z \in \mathcal{N}^\delta} Q^\infty(x_{\text{stable}}, z) + \frac{\eta}{2} \leq -\inf_{z \in N} Q^\infty(x_{\text{stable}}, z) + \eta$$

if $\delta > 0$ is small enough. Therefore, it is sufficient to link the result about the limiting diffusion to the corresponding statement dealing with $\xi^{s,y}$. By Lemma 4.6, we can find $T_1 > 0$, $\varepsilon_1 > 0$ and $r_1 > 0$ such that

$$
\begin{aligned}
(4.16) \quad & \varepsilon \log \sup_{y \in S_{2\varrho}, s \geq r} \mathbb{P}_{s,y}(\Sigma_\varrho > T_1) \\
& \leq -\inf_{z \in N} Q^\infty(x_{\text{stable}}, z) + \frac{\eta}{2} \qquad \forall \varepsilon \leq \varepsilon_1, r \geq r_1.
\end{aligned}
$$

If $\Sigma_\varrho \leq T_1$ and $\rho_{0,T_1}(\xi^{s,y}, Y^\infty) \leq \delta$, then $\{\xi^{s,y}_{\Sigma_\varrho} \in N\}$ is contained in $\{Y^\infty_{\Sigma^\delta_\varrho} \in \mathcal{N}^\delta\}$. Thus,

$$
\begin{aligned}
\mathbb{P}(\xi^{s,y}_{\Sigma_\varrho} \in N) &\leq \mathbb{P}(\xi^{s,y}_{\Sigma_\varrho} \in N, \Sigma_\varrho < T_1) + \mathbb{P}_{s,y}(\Sigma_\varrho \geq T_1) \\
&\leq \mathbb{P}(Y^{\infty,y}_{\Sigma^\delta_\varrho} \in \mathcal{N}^\delta) + \mathbb{P}(\rho_{0,T_1}(\xi^{s,y}, Y^{\infty,y}) \geq \delta) \\
&\quad + \mathbb{P}_{s,y}(\Sigma_\varrho \geq T_1).
\end{aligned}
$$

By (4.16) and Proposition 3.10, the logarithmic asymptotics of the sum on the r.h.s. is dominated by the first term, that is, the lemma is established. □

We are now in a position to establish the lower bound for the exit time which complements Proposition 4.7 and completes the proof of Theorem 4.2.

PROPOSITION 4.9. *There exists $\eta_0 > 0$ such that for any $\eta \leq \eta_0$*

$$
(4.17) \qquad \limsup_{\varepsilon \to 0} \varepsilon \log \mathbb{P}_{x_0}[\tau^\varepsilon_D < e^{(\overline{Q}_\infty - \eta)/\varepsilon}] \leq -\eta/2
$$

*and*

$$
(4.18) \qquad \liminf_{\varepsilon \to 0} \varepsilon \log \mathbb{E}_{x_0}[\tau^\varepsilon_D] \geq \overline{Q}_\infty.
$$

PROOF. In a first step we apply Lemma 4.8 and an adaptation of Lemma 5.7.23 in [7]. The latter explains that the behavior of an Itô diffusion on small time intervals is similar to the behavior of the martingale part, which in our situation is simply given by $\sqrt{\varepsilon} W_t$. We find $r_0 > 0$, $T > 0$ and $\varepsilon_0 > 0$ such that for $\varepsilon \leq \varepsilon_0$

$$
\sup_{y \in S_{2\varrho}, s \geq r_0} \mathbb{P}(\xi^{s,y}_{\Sigma_\varrho} \in \partial D) \leq e^{-(\overline{Q}_\infty - \eta/2)/\varepsilon},
$$

$$
\sup_{y \in D, s \geq r_0} \mathbb{P}\left( \sup_{0 \leq t \leq T} \|\xi^{s,y}_t - y\| \geq \varrho \right) \leq e^{-(\overline{Q}_\infty - \eta/2)/\varepsilon}.
$$

In the sequel, we shall proceed as follows. First, we wait for a large period of time $r_1$ until the diffusion becomes "sufficiently homogeneous," which is possible thanks to the stability assumption. Since $x_{\text{stable}}$ attracts all solutions of the deterministic system, we may find $r_1 \geq r_0$ such that $\psi_r(x_0) \in B_\varrho(x_{\text{stable}})$ for $r \geq r_1$. Second, after time $r_1$, we employ the usual arguments for homogeneous diffusions. Following [7], we recursively define two sequences of stopping times that shall serve to track the diffusion's excursions between the small ball $B_\varrho(x_{\text{stable}})$ around the equilibrium point and the larger sphere $S_{2\varrho} = \partial B_{2\varrho}(x_{\mathrm{stable}})$, before it finally exits from the domain $D$.

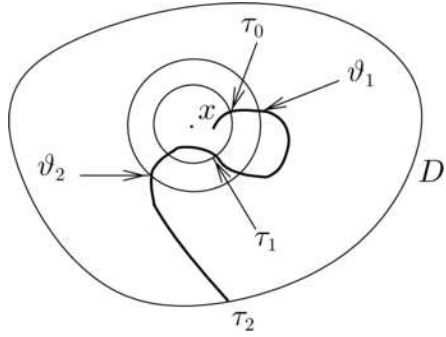

Set $\vartheta_0 = r_1$, and for $m \geq 0$ let

$$\tau_m = \inf\{t \geq \vartheta_m : X_t^\varepsilon \in B_\varrho \cup \partial D\}$$

and

$$\vartheta_{m+1} = \inf\{t > \tau_m : X_t^\varepsilon \in S_{2\varrho}\}.$$

Let us decompose the event $\{\tau_D^\varepsilon \leq kT + r_1\}$. We have

$$\begin{aligned}
(4.19) \qquad \mathbb{P}_{x_0}(\tau_D^\varepsilon \leq kT + r_1) &\leq \mathbb{P}_{x_0}(\{\tau_D^\varepsilon \leq r_1\} \cup \{X_{r_1}^\varepsilon \notin B_{2\varrho}(x_{\mathrm{stable}})\}) \\
&\quad + \sup_{y \in S_{2\varrho}, s \geq r_1} \mathbb{P}_{s,y}(\tau_D^\varepsilon \leq kT).
\end{aligned}$$

The first probability on the r.h.s. of this inequality tends to 0 as $\varepsilon \to 0$. Indeed, by the large deviations principle for $X^\varepsilon$ on the time interval $[0, r_1]$, there exist $\eta_0 > 0$ and $\varepsilon_2 > 0$ such that

$$\varepsilon \log \mathbb{P}_{x_0}(\{\tau_D^\varepsilon \leq r_1\} \cup \{X_{r_1}^\varepsilon \notin B_{2\varrho}(x_{\mathrm{stable}})\}) \leq -\eta/2$$

for $\varepsilon \leq \varepsilon_2$ and $\eta \leq \eta_0$. For the second term in (4.19), we can observe two different cases: either the diffusion exits from $D$ during the first $k$ exits from $D \setminus B_\varrho(x_{\mathrm{stable}})$, or the minimal time spent between two consecutive exits is smaller than $T$. This reasoning leads to the bound

$$\mathbb{P}_{s,y}(\tau_D^\varepsilon \leq kT) \leq \sum_{m=0}^{k} \mathbb{P}_{s,y}(\tau_D^\varepsilon = \tau_m) + \mathbb{P}_{s,y}\left(\min_{1 \leq m \leq k}(\vartheta_m - \tau_{m-1}) \leq T\right).$$

Let us now link these events to the probabilities presented at the beginning of the proof. We have

$$\sup_{y \in S_{2\varrho}, s \geq r_1} \mathbb{P}_{s,y}(\tau_D^\varepsilon = \tau_m) \leq \sup_{y \in S_{2\varrho}, s \geq r_0} \mathbb{P}_{s,y}(\xi_{\Sigma_\varrho}^{s,y} \in \partial D)$$

and

$$\sup_{y \in S_{2\varrho}, s \geq r_1} \mathbb{P}_{s,y}((\vartheta_m - \tau_{m-1}) \leq T) \leq \sup_{y \in S_{2\varrho}, s \geq r_0} \mathbb{P}_{s,y}\left(\sup_{0 \leq t \leq T} \|\xi_t^{s,y} - y\| \geq \varrho\right),$$

which yields the bound

$$\sup_{y \in S_{2\varrho}, s \geq r_1} \mathbb{P}_{s,y}(\tau_D^\varepsilon \leq kT) \leq (2k+1)e^{-(\overline{Q}_\infty - \eta/2)/\varepsilon}.$$

Thus, by choosing $k = \lfloor (e^{(\overline{Q}_\infty - \eta)/\varepsilon} - r_1)/T \rfloor + 1$, we obtain from (4.19)

$$\mathbb{P}_{x_0}(\tau_D^\varepsilon \leq e^{(\overline{Q}_\infty - \eta)/\varepsilon}) \leq e^{-\eta/2\varepsilon} + 5T^{-1}e^{-\eta/2\varepsilon},$$

that is, (4.17) holds. Moreover, by using Chebyshev's inequality, we obtain the claimed lower bound for the expected exit time. Indeed, we have

$$\mathbb{E}_{x_0}(\tau_D^\varepsilon) \geq e^{(\overline{Q}_\infty - \eta)/\varepsilon}(1 - \mathbb{P}_{x_0}(\tau_D^\varepsilon \leq e^{(\overline{Q}_\infty - \eta)/\varepsilon}))$$

$$\geq e^{(\overline{Q}_\infty - \eta)/\varepsilon}(1 - (1 + 5T^{-1})e^{-\eta/2\varepsilon}),$$

which establishes (4.18). □

We end this section with the proof of Theorem 4.3 about the exit location.

PROOF OF THEOREM 4.3.    We use arguments similar to the ones of the preceding proof. Let

$$\overline{Q}_\infty(N) = \inf_{z \in N} Q^\infty(x_{\text{stable}}, z),$$

and assume w.l.o.g. that $\overline{Q}_\infty(N) < \infty$. Otherwise, we may replace $\overline{Q}_\infty(N)$ in the following by some constant larger than $\overline{Q}_\infty$. As in the preceding proof, we may choose $T > 0$, $r_0 > 0$ and $\varepsilon_0 > 0$ such that

$$\sup_{y \in S_{2\varrho}, s \geq r_0} \mathbb{P}_{s,y}(\xi_{\Sigma_\varrho}^{s,y} \in \partial N) \leq e^{-(\overline{Q}_\infty(N) - \eta/2)/\varepsilon} \qquad \forall \varepsilon \leq \varepsilon_0,$$

$$\sup_{y \in D, s \geq r_0} \mathbb{P}_{s,y}\left(\sup_{0 \leq t \leq T} \|\xi_t^{s,y} - y\| \geq \varrho\right) \leq e^{-(\overline{Q}_\infty(N) - \eta/2)/\varepsilon} \qquad \forall \varepsilon \leq \varepsilon_0.$$

It suffices to study the event $A = \{\tau_D^\varepsilon \leq kT + r_0\} \cap \{X_{\tau_D^\varepsilon}^\varepsilon \in N\}$ for positive integers $k$. We see that

$$\mathbb{P}_{x_0}(A) \leq \mathbb{P}_{x_0}(X_{r_0}^\varepsilon \notin B_{2\varrho}(x_{\text{stable}})) + \sup_{y \in S_{2\varrho}, s \geq r_0} \mathbb{P}_{s,y}(\tau_D^\varepsilon \leq kT)$$

$$\leq \mathbb{P}_{x_0}(X_{r_0}^\varepsilon \notin B_{2\varrho}(x_{\text{stable}})) + \sum_{m=0}^{k} \mathbb{P}_{s,y}(\tau_D^\varepsilon = \tau_m, \xi_{\tau_D^\varepsilon}^{s,y} \in N)$$

$$+ \mathbb{P}_{s,y}\left(\min_{1 \leq m \leq k}(\vartheta_m - \tau_{m-1}) \leq T\right)$$

$$\leq \mathbb{P}_{x_0}(X_{r_0}^\varepsilon \notin B_{2\varrho}(x_{\text{stable}})) + (2k+1)e^{-(\overline{Q}_\infty(N) - \eta/2)/\varepsilon}.$$

The choice $k = \lfloor (e^{(\overline{Q}_\infty(N) - \eta)/\varepsilon} - r_0)/T \rfloor + 1$ yields

$$\mathbb{P}_{x_0}(A) \le \mathbb{P}_{x_0}(X_{r_0}^\varepsilon \notin B_{2\varrho}(x_{\text{stable}})) + 5T^{-1} e^{-\eta/2\varepsilon}.$$

This implies that $\mathbb{P}_{x_0}(\tau_D^\varepsilon \le e^{(\overline{Q}_\infty(N) - \eta)/\varepsilon}, X_{\tau_D^\varepsilon}^\varepsilon \in N) \to 0$ as $\varepsilon \to 0$. Now choose $\eta$ small enough such that $\overline{Q}_\infty(N) - \eta > \overline{Q}_\infty + \eta$. Then Proposition 4.7 states that the exit time of the domain $D$ is smaller than $e^{(\overline{Q}_\infty + \eta)/\varepsilon}$ with probability close to 1. The combination of these two results implies $\mathbb{P}_{x_0}(X_{\tau_D^\varepsilon}^\varepsilon \in N) \to 0$ as $\varepsilon \to 0$. $\quad\square$

## 5. The gradient case: examples.

The structural assumption about $\Phi$, namely its rotational invariance as stated in (2.4), implies that $\Phi$ is always a potential gradient. In fact, this assumption means that $\Phi$ is the gradient of the positive potential

$$\mathcal{A}(x) = \int_0^{\|x\|} \phi(u) \, du.$$

In this section, we make the additional assumption that the second drift component given by the vector field $V$ is also a potential gradient, which brings us back to the very classical situation of gradient type time-homogeneous Itô diffusions. In this situation, quasi-potentials and exponential exit rates may be computed rather explicitly and allow for a good illustration of the effect of self-stabilization on the asymptotics of exit times.

We assume from now on that $V = -\nabla U$ is the gradient of a potential $U$ on $\mathbb{R}^d$. Then the drift of the limiting diffusion $Y^\infty$ defined by (3.13) is also a potential gradient, that is,

$$b(x) := V(x) - \Phi(x - x_{\text{stable}}) = -\nabla(U(x) + \mathcal{A}(x - x_{\text{stable}})).$$

A simple consequence of Theorem 3.1 in [8] allows one to compute the quasi-potential explicitly in this setting.

LEMMA 5.1. *Assume that $V = -\nabla U$. Then for any $z \in D$,*

$$Q^\infty(x_{\text{stable}}, z) = 2(U(z) - U(x_{\text{stable}}) + \mathcal{A}(z - x_{\text{stable}})).$$

*In particular,*

$$\overline{Q}_\infty = \inf_{z \in \partial D} 2(U(z) - U(x_{\text{stable}}) + \mathcal{A}(z - x_{\text{stable}})).$$

Observe that the exit time for the self-stabilizing diffusion is strictly larger than that of the classical diffusion defined by

$$dZ_t^\varepsilon = V(Z_t^\varepsilon) \, dt + \sqrt{\varepsilon} \, dW_t, \qquad Z_0^\varepsilon = x_0.$$

Indeed, by the theory of Freidlin and Wentzell,

$$\lim_{\varepsilon \to 0} \varepsilon \log \mathbb{E}_{x_0}(\tau_D^\varepsilon(Z^\varepsilon)) = \inf_{z \in \partial D} 2(U(z) - U(x_{\text{stable}})) \tag{5.1}$$

$$< \overline{Q}_\infty = \lim_{\varepsilon \to 0} \varepsilon \log \mathbb{E}_{x_0}(\tau_D^\varepsilon(X^\varepsilon)). \tag{5.2}$$

The exit problem is in fact completely different if we compare the diffusions with and without self-attraction. We have already seen that the exponential rate is larger in the attraction case. Let us next see by some examples that the exit location may change due to self-stabilization.

5.1. *The general one-dimensional case.* In this subsection we confine ourselves to one-dimensional self-stabilizing diffusions. In dimension one, the structural assumptions concerning $\Phi$ and $V$ are always granted, and we may study the influence of self-stabilization on exit laws in a general setting. Let $a < 0 < b$, and assume for simplicity that the unique stable equilibrium point is the origin 0. Denote by $U(x) = -\int_0^x V(u)\,du$ the potential that induces the drift $V$. As seen before, the interaction drift is the gradient of the potential $\mathcal{A}(x) = \int_0^{|x|} \phi(u)\,du$. Since we are in the gradient situation, the exponential rate for the mean exit time from the interval $[a, b]$ can be computed explicitly.

If we denote by $\tau_x(X^\varepsilon) = \inf\{t \geq 0 \colon X_t^\varepsilon = x\}$ the first passage time of the level $x$ for the process $X^\varepsilon$ and $\tau_I = \tau_a \wedge \tau_b$, then the exit law of the classical diffusion $Z^\varepsilon$ (i.e., without self-stabilization) is described by

$$\lim_{\varepsilon \to 0} \mathbb{P}_0(e^{(Q_0^\infty - \eta)/\varepsilon} < \tau_I(Z^\varepsilon) < e^{(Q_0^\infty + \eta)/\varepsilon}) = 1$$

and

$$\lim_{\varepsilon \to 0} \varepsilon \log \mathbb{E}_0(\tau_I(Z^\varepsilon)) = Q_0^\infty,$$

where $Q_0^\infty = 2\min(U(a), U(b))$. Moreover, if we assume that $U(a) < U(b)$, then $\mathbb{P}_0(\tau_I(Z^\varepsilon) = \tau_a(Z^\varepsilon)) \to 1$ as $\varepsilon \to 0$.

The picture changes completely if we introduce self-stabilization. The quasi-potential becomes

$$Q_1^\infty = 2\min(U(a) + \mathcal{A}(a), U(b) + \mathcal{A}(b)) > Q_0^\infty,$$

so the mean exit time of $X^\varepsilon$ from the interval $I$ is strictly larger compared to that of $Z^\varepsilon$. This result corresponds to what intuition suggests: the process needs more work and consequently more time to exit from a domain if it is attracted by some law concentrated around the stable equilibrium point. Furthermore, if $a$ and $b$ satisfy

$$\mathcal{A}(b) - \mathcal{A}(a) < U(a) - U(b),$$

we observe that $\mathbb{P}(\tau_I(X^\varepsilon) = \tau_b(X^\varepsilon)) \to 1$, that is, the diffusion exits the interval at the point $b$. Thus, we observe the somehow surprising behavior that self-stabilization changes the exit location from the left to the right endpoint of the interval. See Figure 1.

### 5.2. An example in the plane.
In this subsection, we give another explicit example in dimension two, in order to illustrate changes of exit locations in more detail.

Let $V = -\nabla U$, where

$$U(x,y) = 6x^2 + \tfrac{1}{2}y^2,$$

and let us examine the exit problem for the elliptic domain

$$D = \{(x,y) \in \mathbb{R}^2 : x^2 + \tfrac{1}{4}y^2 \le 1\}.$$

The unique stable equilibrium point is the origin $x_{\text{stable}} = 0$.

The asymptotic mean exit time of the diffusion $Z_t^\varepsilon$ starting in 0 is given by $\lim_{\varepsilon \to 0} \varepsilon \log \mathbb{E}_0(\tau_{\bar{D}}^\varepsilon(Z^\varepsilon)) = 4$, since the minimum of the potential on $\partial D$ is reached if $y = \pm 2$ and $x = 0$. Let us now focus on its exit location, and denote $N_{(x,y)} = \partial D \cap B_\varrho((x,y))$. The diffusion exits asymptotically in the neighborhood $N_{(0,2)}$ with probability close to $1/2$ and in the neighborhood $N_{(0,-2)}$ with the same probability.

Now we look how self-stabilization changes the picture. For the interaction drift we choose $\Phi(x,y) = \nabla \mathcal{A}(x,y)$, with $\mathcal{A}(x,y) = 2x^2 + 2y^2$. First, the self-stabilizing diffusion $X^\varepsilon$ starting in 0 needs more time to exit from $D$, namely $\lim_{\varepsilon \to 0} \varepsilon \log \mathbb{E}_0(\tau_{\bar{D}}^\varepsilon(X^\varepsilon)) = 16$. More surprisingly, though, the exit location is completely different. The diffusion exits asymptotically with probability close to $1/2$ in the neighborhoods $N_{(-1,0)}$ and $N_{(1,0)}$, respectively.

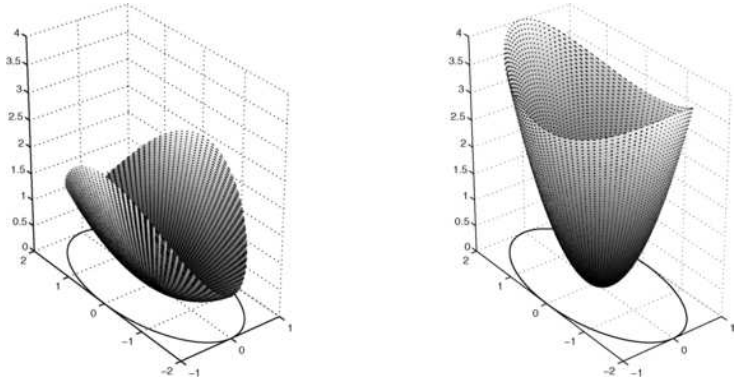

Fig. 1. *Potentials $U$ (left picture) and $U + \mathcal{A}$ (right picture).*

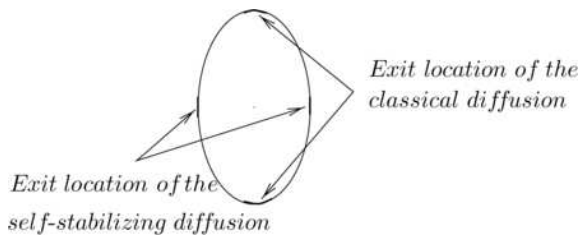

*Exit location of the classical diffusion*

*Exit location of the self-stabilizing diffusion*

**Acknowledgment.** We are very much indebted to an anonymous referee for his correct reading and his useful remarks.

Ecole des Mines de Nancy
and Institut de Mathématiques Elie Cartan
BP 239
54506 Vandoeuvre-lès-Nancy
France
E-mail: herrmann@iecn.u-nancy.fr

Institut für Mathematik
Humboldt-Universität zu Berlin
Unter den Linden 6
10099 Berlin
Germany
E-mail: imkeller@mathematik.hu-berlin.de
        peithman@mathematik.hu-berlin.de